\documentclass[a4paper,10pt]{article}
\usepackage[left=3cm,top=2cm,bottom=2cm,right=2.2cm,includehead,includefoot]{geometry}
\usepackage[ansinew]{inputenc}
\usepackage[american]{babel}
\usepackage{amsfonts} 
\usepackage{amssymb}
\usepackage{amsthm}
\usepackage{amsmath}
\usepackage{mathtools}
\usepackage[arrow, matrix, curve]{xy}
\usepackage{enumitem}
\usepackage{dsfont}
\usepackage{bbm}
\usepackage{url}

\usepackage{mathrsfs}

\newcommand{\ov}[1]{\overline{#1}}
\newcommand{\un}[1]{\underline{#1}}
\newcommand{\lmat}{\left(\begin{array}}
\newcommand{\rmat}{\end{array}\right)}
\newcommand{\M }[1]{\mathbb #1}
\newcommand{\MC}[1]{\mathcal{#1}}

\newcommand{\unif}{\mathcal{\pi}}
\newcommand{\Sp}{{\rm Spec} \,}

\newcommand{\op}[1]{\operatorname{#1}}
\newcommand{\exit}{\hfill $\square$}

\newtheorem{lemma}{Lemma}[section]
\newtheorem{theorem}[lemma]{Theorem}
\newtheorem{definition}[lemma]{Definition}
\newtheorem{corollary}[lemma]{Corollary}
\newtheorem{proposition}[lemma]{Proposition}
\newtheorem{conjecture}[lemma]{Conjecture}
\newtheorem{construction}[lemma]{Construction}
\newtheorem{notation}[lemma]{Notation}

\newcommand{\thm}[3]{\begin{theorem} #1 #2 $\left.\right.$ \\ #3 \end{theorem}}
\newcommand{\prop}[3]{\begin{proposition} #1 #2 $\left.\right.$ \\ #3 \end{proposition}}
\newcommand{\lem}[3]{\begin{lemma} #1 #2 $\left.\right.$ \\ #3 \end{lemma}}
\newcommand{\cor}[3]{\begin{corollary} #1 #2 $\left.\right.$ \\ #3 \end{corollary}}
\newcommand{\defi}[3]{\begin{definition} #1 #2 $\left.\right.$ \\ #3 \end{definition}}
\newcommand{\conj}[3]{\begin{conjecture} #1 #2 $\left.\right.$ \\ #3 \end{conjecture}}

\newtheorem{preremark}[lemma]{Remark}  
\newenvironment{remark}{\begin{preremark}\rm}{\end{preremark}}
\newcommand{\rem}[2]{\begin{remark} #1 $\left.\right.$ \\ #2 \end{remark}}

\newcommand{\prooof}{\proof $\left. \right.$ \\}

\usepackage{hyperref}

\title{A refinement of the Hodge stratification \\ for connected reductive groups}
\author{Stephan Neupert}

\begin{document}
\maketitle

\begin{abstract}
  For connected reductive groups $G$ over a finite extension $F$ of $\M{Q}_p$ and $L$ the maximal unramified extension of $F$ we study the sets $H_{\un{\mu}, N}(G)$ of elements $b \in G(L)$ with given Hodge points of $(b\sigma), (b\sigma)^2, \ldots, (b\sigma)^N$. We explain the relationship to stratifications of some moduli scheme of abelian varieties defined by Goren and Oort respectively Andreatta and Goren. We show that for sufficiently large $N$ the Newton point is constant on the sets $H_{\un{\mu}, N}(G)$ and compute such $N$ for certain classes of groups.
\end{abstract}

\section{Introduction}\label{sec:Introduction}
Let $\MC{M}$ be a moduli scheme over $\M{F}_p$ of abelian varieties of PEL-type. To get a better understanding of the geometry of such moduli schemes one aims for a refinement of the stratification of $\MC{M}$ according to the $p$-rank of the abelian variety's $p$-torsion. \\
One way to do this is the stratification via Newton points. This invariant of any abelian variety over an algebraically closed field $k$ of positive characteristic is constant on isogeny classes and can be constructed in the following way: To an abelian variety $A$ of dimension $g$ associate its $p$-divisible group $A[p^\infty] = \varinjlim A[p^n]$ where $A[p^n]$ denotes the group of $p^n$-torsion points in $A$. It is a limit of finite locally free group schemes over $k$. The category of such group schemes is equivalent to the opposite category of Dieudonn\'e modules. These are free modules $D$ over the Witt-ring $W(k)$ over $k$, endowed with a Frobenius-linear morphism $F$ and a Frobenius-antilinear morphism $V$ satisfying $FV = VF = p$ (where $p$ denotes the multiplication with $p \in W(k)$). Note that the PEL-data on the abelian variety give corresponding structures on the Dieudonn\'e module. Denoting the Frobenius on $k$ by $\sigma$, the morphism $F$ can be identified with a $\sigma$-conjugacy class in the reductive subgroup of $GL_g$ of morphisms respecting these additional structures. Kottwitz gave a purely group-theoretic description of the $\sigma$-conjugacy classes extending the notion of classical Newton points for $GL_g$ (which classify $F$ after forgetting about the PEL-structure) to all connected reductive groups (cf. \cite{KottIso1} and \cite{KottIso2}). \\
Another refinement of the $p$-rank stratification was introduced by Goren and Oort in \cite{GorenOortStrata} (for the unramified case) respectively by Andreatta and Goren in \cite{AndreattaGoren} (for the totally ramified case). In this introduction let us consider only the second situation: Let $\MC{M}(\M{F}_p, \mu_N)$ be the fine moduli scheme over $\M{F}_p$ of polarized abelian varieties with real multiplication by the ring of integers $O_{\tilde{F}}$ of a totally real extension $\tilde{F}/\M{Q}$ of degree $g$ and totally ramified at $p$ and a $\mu_N$-level structure (with $N \geq 4$) satisfying the Deligne-Pappas condition. This forces those abelian varieties to have dimension exactly $g$. Andreatta and Goren notice that, contrary to the unramified case, the Ekedahl-Oort-stratification (using some discrete invariant of the $p$-torsion of abelian varieties, cf. \cite{OortStrata}) does not even define a stratification in this situation. Instead they introduce two invariants: To get a definition allowing explicit computations, associate to the abelian variety $A$ over $k$ its display $(P, Q, V^{-1}, F)$. Then one can define a normal form for such displays in which one can find both invariants $j$ and $n$ as exponents (cf. \cite[\S4.4]{AndreattaGoren} or loc. cit. \S \ref{subsec:AndreattaGoren}). Their motivation to define $j$ and $n$ comes from the following description, which is pretty similar to the way the invariants in the unramified setting are defined in \cite{GorenOortStrata}: The first deRham-cohomology $H^1_{{\rm dR}}(A)$ is a free $k[T]/(T^g)$-module of rank $2$. Then we can find generators $\alpha$ and $\beta$ of $H^1_{{\rm dR}}(A)$ such that:
\[H^1(A, \MC{O}_A) = (T^i\alpha) \oplus (T^j\beta) \quad, \quad i \geq j, i + j = g\] 
(viewing $H^1(A, \MC{O}_A)$ as a quotient of $H^1_{{\rm dR}}(A)$). This gives the first invariant $j$. For the second let the $a$-number $a(A)$ be the nullity of the Hasse-Witt matrix of $A$. Then put $n = a(A) - j$. Now decompose $\MC{M}(\M{F}_p, \mu_N)$ into the loci $\MC{M}_{j, n}$ (with $0 \leq j \leq n \leq g$) of abelian varieties with these invariants (cf. \cite[\S5]{AndreattaGoren}). Andreatta and Goren now prove that these $\MC{M}_{j, n}$ indeed define a stratification with many nice properties, e.g. that $j$ and $n$ determine the Newton point corresponding to the abelian varieties in $\MC{M}_{j, n}$ and that the action of the Hecke correspondence on the strata can be described explicitly. \\
The aim of this work is to give a purely group-theoretic description of both the stratification defined by Goren and Oort and the one of Andreatta and Goren and to extend it to all connected reductive groups $G$ in such a way that it has similar properties with respect to Newton points. \\
To define this partition for a reductive group $G$ over $O_F$ (the ring of integers of a finite field extension $F$ of $\M{Q}_p$) with maximal torus $T$ and Weyl group $\Omega$ let $L$ be the maximal unramified extension of $F$ with ring of integers $O_L$ and define $\mu: G(O_L)\backslash G(L)/G(O_L) \to X_*(T)/\Omega$ to be the map associating to each element its Hodge point (cf. \S\ref{subsec:DefiHodgeNewton}). Fix a positive integer $N$ and a $N$-tuple $\un{\mu} = (\mu^1, \mu^2, \ldots, \mu^N)$ of Hodge points. Then let $H_{\un{\mu}, N}(G) \subset G(L)$ be the set of elements $b \in G(L)$ with $\mu((b\sigma)^i) = \mu^i$ for $i = 1, \ldots, N$ (again denoting the Frobenius by $\sigma$).\\
The sections \ref{sec:GLn} and \ref{sec:GeneralG} discuss the connection between the partition $H_{\un{\mu}, k}(G)$ and Newton points first for the case $G = GL_n$ and then for arbitrary connected reductive groups. The main result here states that there is a constant $C_{\mu, G}$ depending on a Hodge point $\mu$ (and of course on $G$) such that the Newton point is constant on each $H_{\un{\mu}, C_{\mu, G}}(G)$ where $\un{\mu}$ is a $C_{\mu, G}$-tuple with first entry $\mu$. To show this we analyze the convergence of the sequence of Hodge points $(\frac 1i \mu((b\sigma)^i))_i$ to the Newton point $\nu(b)$. \\
The partition of $G(L)$ in the $H_{\un{\mu}, k}(G)$ is most useful in the case of small constants $C_{\mu, G}$, in particular when one can take $C_{\mu, G} = 2$. In section \ref{sec:Explicit} we give some explicit constants $C_{\mu, G}$ for special classes of groups. In particular we obtain a result which turns out to be a slight generalization of Andreatta and Goren's theorem that the Newton point is constant on the strata $\MC{M}_{j, n}$. \\
Finally we study the relationship with the stratifications mentioned above: We see in section \ref{subsec:GorenOort} that the strata defined by Goren and Oort in \cite{GorenOortStrata} can be described as the loci where the Frobenius of the Dieudonn\'e module associated to the abelian variety lies in $H_{\un{\mu}, 2}({\rm Res}_{F/\M{Q}_p}(GL_2))$ for suitable $\un{\mu}$. In section \ref{subsec:AndreattaGoren} the same is done in the totally ramified case.

\subsection*{Acknowledgements}
I thank Eva Viehmann for her suggestion to study this stratification and many helpful discussions. I also thank the referee for all the comments, especially for mentioning the possible relationship to \cite{GorenOortStrata}. \\
The author was partially supported by ERC starting grant 277889 'Moduli spaces of local $G$-shtukas'.

\section{Definition and first properties}\label{sec:Definitions}
\subsection{Setting}\label{subsec:Setting}
Throughout the whole paper we will fix the following:
\begin{itemize}
\item $k$ an algebraically closed field of characteristic $p > 0$
\item $K = {\rm Frac}(W(k))$ the fraction field of the ring of Witt vectors over $k$
\item $F$ a finite extension of $\M{Q}_p$ inside a fixed algebraic closure $\ov{K}$ of $K$ 
\item $L = K \cdot F$ the composite of $K$ and $F$ inside $\ov{K}$
\item $\sigma$ the Frobenius automorphism of $L/F$ induced by the Frobenius on their residue fields
\item $\Gamma = Gal(\ov{F}, F)$
\item $G$ a connected reductive linear algebraic group over $F$
\item $T \subset G$ a maximal torus containing a maximal $L$-split torus
\item $\Omega$ the Weyl group associated to $T$ 
\item $B(G)$ the set of $\sigma$-conjugacy classes of $G(L)$
\item $\Lambda_0$ a special vertex in the Bruhat-Tits building of $G$ over $L$
\end{itemize}

\rem{}{
i) All valuations of fields are assumed to have value group $\M{Z}$. \\
ii) A similar theory with analogous proofs should work in the function field case, i.e. $K = k((t))$. \\
iii) In section \ref{sec:OtherStrata} $F$ will also denote the Frobenius morphism of a Dieudonn\'e module. We hope that it will always be clear what is meant from the context.
}

\subsection{Hodge and Newton points}\label{subsec:DefiHodgeNewton}
We will recall the definition of the Hodge and Newton points for arbitrary connected reductive groups $G$ over $F$. \\
For the Newton point we state the explicit description given in Kottwitz \cite[\S4.3]{KottIso1}: \\
Let $\M{D}$ be the pro-algebraic torus over $\M{Q}_p$ with character group $\M{Q}$. For $b \in G(L)$ let $\nu \in {\rm Hom}_L(\M{D}, G)$ be the unique element for which there is an integer $n > 0$, a uniformizer $\unif \in L$ and an element $c \in G(L)$ such that the following three conditions hold: 

i) $n\nu \in {\rm Hom}_L(\M{G}_m, G)$

ii) ${\rm Int}(c) \circ (n\nu)$ is defined over the fixed field of $\sigma^n$ on $L$ 

iii) $c (b\sigma)^n c^{-1} = c \cdot (n\nu)(\unif) \cdot c^{-1} \cdot \sigma^n$ \\
(where ${\rm Int}(c)$ denotes the conjugation by $c$). This defines a map $\nu: G(L) \to {\rm Hom}_L(\M{D}, G)$. Two $\sigma$-conjugate elements in $G(L)$ are mapped to elements which are conjugate under $G(L)$. Furthermore the elements in the image of the map $\nu$ are invariant under the action of the Frobenius $\sigma$. Hence this induces a map defining the Newton point
\[\nu: B(G) \to ({\rm Int} \; G(L) \backslash {\rm Hom}_L(\M{D}, G))^{\sigma}.\]
As all maximal tori of $G$ are conjugate, fixing one of them, namely $T$, gives an isomorphism 
\[\left({\rm Int} \; G(L) \backslash {\rm Hom}_L(\M{D}, G)\right)^{\sigma} \cong (X_*(T)_{\M{Q}} / \Omega)^\Gamma.\]
This allows us to view $\nu$ as a map
\[\nu: B(G) \to (X_*(T)_{\M{Q}} / \Omega)^\Gamma.\]

\rem{\label{rem:GLNormalForm}}{
If $G = GL_n$ then the conditions i), ii) and iii) give the classical Dieudonn\'e-Manin classification of isocrystals after choosing some basis and representing the isocrystal by a matrix: For $m \in \M{N}$, $h \in \M{Z}$ with $gcd(m, h) = 1$, let $B_{m, h} \in GL_m$ with $B_{m, h}(e_i) = e_{i + 1}$ for $i = 1,\ldots, m-1$ and $B_{m, h}(e_m) = \unif^h e_1$ (where $\{e_i\}_{i = 1, \ldots, m}$ are basis vectors). Then any $b \in GL_n(L)$ can be $\sigma$-conjugated to a block-matrix $b_0$ with blocks of the form $B_{m, h}$. In the following $b_0$ is called the normal form of $b \in GL_n$ for the Newton point. This block-matrix $b_0$ is unique up to permutation of the blocks and each block $B_{m, h}$ gives $m$ slopes $\frac hm$. \\
To compare this to the general Newton point assume that $T$ is the diagonal torus in $GL_n$. Then $\nu(b)$ is the (unique up to $\Omega = S_n$-action) rational cocharacter such that $(n!\nu(b))(\unif) = b_0^{n!}$.
}
$\left.\right.$ \vspace{1mm} \\
For the Hodge point we follow \cite[\S3.3]{Tits}: \\
Let $G(O_L) \subset G(L)$ be the stabilizer of the fixed special vertex $\Lambda_0$ in the Bruhat-Tits building (in analogy to the classical case $G = GL_n$). Let $S$ be a maximal $L$-split torus of $G$ contained in $T$ (this exists by our choice of $T$) and denote the centralizer of $S$ by $Z$. Denoting the valuation on $L$ by ${\rm val}_L$ define $v: Z(L) \to X_*(S)_{\M{Q}} \subseteq X_*(T)_{\M{Q}}$ via the equality
\[\langle \chi, v(z) \rangle = -{\rm val}_L(\chi(z)) \qquad {\rm for \; all \;} z \in Z(L), \chi \in X^*(Z) \]
(as in the definition of an apartment in the Bruhat-Tits building, cf. \cite[\S 1.2]{Tits}). Then this map induces a well-defined map (cf. \cite[\S4]{BT}):
\[\begin{array}{cccc}
    \mu \coloneqq \mu_{\Lambda_0}: & G(O_L) \backslash G(L) / G(O_L) & \longrightarrow & X_*(T)_\M{Q}/\Omega \\
      & G(O_L)zG(O_L) & \; \xmapsto{\quad}{} & -v(z)
  \end{array}\]
where $z \in Z(L)$. \\
For an alternative way to describe $\mu$ for unramified reductive groups, cf. \cite[\S4]{RapoRich} and \cite{KottSVTOI}.

\rem{}{
i) The Hodge point depends on the choice of the special vertex. If $\Lambda_0' = g \Lambda_0$ with $g \in G(L)$ is another such vertex, then the Hodge points compare via
\[\mu_{\Lambda_0'}(b) = \mu_{\Lambda_0}(g^{-1}b\sigma(g)).\] 
ii) For $G = GL_n$ identify $G = GL(V)$ for some $F$-vector space $V$. Then the $G(L)$-orbit $\MC{V}$ of $\Lambda_0$ in the Bruhat-Tits building consists of the $O_L$-lattices in $V$. Thus the stabilizer of $\Lambda_0 \in \MC{V}$ can be identified with the set of matrices $GL_n(O_L)$ invertible over $O_L$ (after fixing a basis of the lattice $\Lambda_0$). Furthermore $T = S = Z$. Then one computes:
\[(-v)({\rm diag}(t_1, \ldots, t_n))(\lambda) = {\rm diag}(\lambda^{{\rm val}_L(t_1)}, \ldots, \lambda^{{\rm val}_L(t_n)}) \qquad {\rm for \; all \;} \lambda \in L^\times\]
Thus the Hodge slopes of an element $b \in G(O_L) \cdot {\rm diag}(t_1, \ldots, t_n) \cdot G(O_L)$ are exactly the valuations of $t_1, \ldots, t_n$. But these are by definition the elementary divisors of $b$. Hence for $G = GL_n$ the construction above recovers the classical definition of Hodge points.
}
$\left.\right.$ \vspace{2mm} \\
Furthermore we define the following partial ordering $\prec$ on $X_*(T)_{\M{R}} / \Omega$ (cf. e.g. \cite[lemma 2.1]{RapoRich}): \\
For a fixed basis of the set of roots in $X^*(T)$ let $\ov{C} \subset X_*(T)_{\M{R}}$ be the Weyl chamber and $C^\vee \subset X_*(T)_{\M{R}}$ be the obtuse Weyl chamber. Then say $x \prec x'$ for $x, x' \in X_*(T)_{\M{R}} / \Omega$ if for the representatives $\tilde{x}$ resp. $\tilde{x}'$ of $x$ resp. $x'$ in $\ov{C} \subset X_*(T)$ one of the following equivalent conditions is satisfied: \\
i) $x$ lies in the convex hull of $\{\omega x'; \omega \in \Omega\}$. \\
ii) $\tilde{x}' - \tilde{x} \in C^\vee$. \\
iii) $\tilde{x}' - \omega \tilde{x} \in C^\vee$ for all $\omega \in \Omega$. \\
iv) For any representation $\rho: G \to GL(V)$ and maximal torus $T' \in GL(V)$ containing $\rho(T)$ we have $\rho(x) \prec \rho(x')$ (using one of the other three conditions or the reformulation in the following remark to define $\prec$ for $GL_n$).

\rem{}{
Consider the special case $G = GL_n$, $B \subset G$ the Borel subgroup of upper triangular matrices, $T \subset B$ the diagonal torus and the usual identification $X_*(T)_{\M{R}} \cong \M{R}^n$ via the coefficients on the diagonal. For the usual choice of simple roots, a representative $\tilde{x} \in X_*(T)_{\M{R}}$ lies in the dominant chamber if and only if its corresponding element $(\lambda_1, \ldots, \lambda_n) \in \M{R}^n$ fulfills
\[\lambda_1 \geq \lambda_2 \geq \ldots \geq \lambda_n\]
Then two elements $x, x' \in X_*(T)_{\M{R}} / \Omega$ satisfy $x \prec x'$ if and only if the corresponding dominant representatives $(\lambda_1, \ldots, \lambda_n)$ resp. $(\lambda_1', \ldots, \lambda_n')$ fulfill:
\[\sum_{i = 1}^h \lambda_1 \leq \sum_{i = 1}^h \lambda_i' \qquad {\rm for} \; h = 1, \ldots, n\]
with equality for $h = n$.
}

\defi{}{}{
Fix a basis $\{\alpha_i\}_{i = 1, \ldots, n} \subset X^*(T)_{\M Q}$. Let $x, x' \in X_*(T)_{\M{Q}} / \Omega$ with representatives $\tilde{x}$ resp. $\tilde{x}'$ in $\ov{C}$. Then define:
\[|x, x'| = \max_{I \subset \{1, \ldots, n\}}\left\{\left|\left\langle \tilde{x} - \tilde{x}', \sum_{\alpha_i \in I} \alpha_i \right\rangle\right|\right\}\]
}

\rem{\label{rem:MetricGLn}}{
i) The definition of $|x, x'|$ depends on the $\alpha_i$. But as all such norms are equivalent it will not matter which basis we actually fix. \\
ii) If $G = GL_n$, $T \subset G$ the diagonal torus and $\alpha_i$ the projection to the $i$th coordinate of $T$. Then $\langle x, \alpha_i \rangle$ is the $i$-th slope of $x$. Assume now that the sum of all slopes for $x$ and $x'$ coincide (this will be the only case of interest for us) and let $(\lambda_1, \ldots, \lambda_n)$ resp. $(\lambda_1', \ldots, \lambda_n')$ be representatives of $x$ resp. $x'$ in $\ov{C}$, i.e. the $\lambda_i$ resp. $\lambda_i'$ are monotonically decreasing. Then
\[|x, x'| = \max_{I \subset \{1, \ldots, n\}} \left\{\left|\sum_{i \in I} (\lambda_i - \lambda_i')\right|\right\} = \sum_{\lambda_i > \lambda_i'} (\lambda_i - \lambda_i')\]
In particular $|x, x'|$ is at least the maximal vertical distance between the corresponding concave polygons for $x$ resp. $x'$.
}

\subsection{Definition and admissibility of \texorpdfstring{$\mathbf{H_{\un{\mu}, k}(G)}$}{refined Hodge strata}} \label{subsec:DefiRefinedHodgeStrata}
Before defining the stratification $H_{\un{\mu}, k}(G)$, recall the classical Hodge and Newton stratifications:
\defi{}{}{
i) The Hodge stratum for $\mu \in X_*(T)_{\M{Q}} / \Omega$ is:
\[H_{\mu}(G) = \{b \in G(L) \; | \; \mu(b) = \mu\}\]
ii) The Newton stratum for $\nu \in X_*(T)_{\M{Q}} / \Omega$ is:
\[M_{\nu}(G) = \{b \in G(L) \; | \; \nu(b) = \nu\}\]
}

\defi{}{}{
Let $k \geq 1$ and $\un{\mu} = (\mu^1, \mu^2, \ldots, \mu^k) \in (X_*(T)_{\M{Q}} / \Omega)^k$. Then let
\[H_{\un{\mu}, k}(G) = \{b \in G(L) \; | \; \mu((b\sigma)^i) = \mu^i \; {\rm for} \; i = 1, 2, \ldots, k \}\]
where $(b\sigma)^i = b \cdot \sigma(b) \cdot \sigma^2(b) \cdot \ldots \cdot \sigma^{i-1}(b)$. We call $H_{\un{\mu}, k}(G)$ a refined Hodge stratum in $G(L)$.
}

\rem{}{
Although we name these sets 'strata', they only form a partition of $G(L)$ but miss the topological properties of actual stratifications.
}

\lem{\label{lem:DefiAdmissible}}{}{
For any two Hodge points $\mu, \mu' \in X_*(T)_{\M{Q}} / \Omega$ there is a $\sigma$-invariant bounded open subgroup $K_{\mu, \mu'} \subset G(L)$ contained in $G(O_L)$ and depending only on $\mu$ and $\mu'$ such that for all $b \in H_{\mu}(G)$, $b' \in H_{\mu'}(G)$ and $g \in K_{\mu, \mu'}$:
\[\mu(b \cdot g \cdot b') = \mu(b \cdot b')\]
}

\prooof
Fix for now any representation $\rho: G(L) \longrightarrow GL_n(L)$ mapping $G(O_L)$ into $GL_n(O_L)$ (the $O_L$-matrices which are invertible over $O_L$).
Then $\rho(H_{\mu'}(G)) \subset H_{\rho(\mu')}(GL_n)$ for some Hodge point $\rho(\mu')$. Choose a positive integer $N$ with $N > \langle\rho(\mu'), \alpha \rangle$ for any root $\alpha$ of $GL_n$ and consider $K_{\rho(\mu), \rho(\mu')} = \{\tilde{g} \in GL_n(O_L) \; | \; \tilde{g} \equiv 1 \bmod \unif^N\}$. Then $K_{\rho(\mu), \rho(\mu')}$ is indeed a $\sigma$-invariant open bounded normal subgroup contained in $GL_n(O_L)$. Furthermore by choice of $N$ for any $\tilde{b}' \in H_{\rho(\mu')}(GL_n)$ and $\tilde{g} \in K_{\rho(\mu), \rho(\mu')}$ we have $\tilde{b}'^{-1} \cdot \tilde{g} \cdot \tilde{b}' \in G(O_L)$. This implies for any $b, b'$ and $g$ as in the statement:
\[b \cdot g \cdot b' = b \cdot b' \cdot (b'^{-1}gb') \in b \cdot b' \cdot G(O_L)\]
and in particular the equality of Hodge points. \\
Now choose a finite set of representations $\{\rho_i\}_i$ such that any two of the finitely many Hodge points appearing as $\mu(bgb')$ for $b \in H_{\mu}(G)$, $b' \in H_{\mu'}(G)$ and $g \in G(O_L)$ can be distinguished by their image under at least one $\rho_i$. The existence of such a set follows from part iv) of the definition of the partial ordering $\prec$ (cf. \S\ref{subsec:DefiHodgeNewton}). Then define:
\[K_{\mu, \mu'} = G(O_L) \cap \bigcap_i \rho_i^{-1}(K_{\rho_i(\mu), \rho_i(\mu')}) \subset G(O_L)\] 
This is an open bounded subgroup contained in $G(O_L)$ such that for any $b \in H_{\mu}(G)$, $b' \in H_{\mu'}(G)$ and $g \in K_{\mu, \mu'}$ the element $bgb'$ and $bb'$ have the same Hodge points after applying any of the $\rho_i$. Hence they have the same Hodge points in $G$. Note that $\sigma(K_{\mu, \mu'}) = K_{\mu, \mu'}$ as it fixes all of the $K_{\rho_i(\mu), \rho_i(\mu')}$.
\exit

\prop{}{}{
Let $G$ be an arbitrary connected reductive group. \\
i) $G(L) = \bigcup_{\un{\mu}} H_{\un{\mu}, k}(G)$ defines a partition of $G(L)$ for fixed $k$. \\
ii) For each $k \in \M{N}$ and $\un{\mu}$ there is a bounded open subgroup $K_{\un{\mu}} \subset G(L)$ with
\[b \cdot K_{\un{\mu}} \subset H_{\un{\mu}, k}(G)\]  
for every $b \in H_{\un{\mu}, k}(G)$.
}

\prooof
i) trivial. \\
ii) For $k = 1$ one may take $K_{\un{\mu}} = G(O_L)$. For $k > 1$ and $\un{\mu} = (\mu^1, \ldots. \mu^k)$ define
\[K_{\un{\mu}} = \bigcap_{\substack{i, j > 0\\i + j \leq k}} K_{\mu^i, \sigma^i(\mu^j)} \subset G(O_L)\]
with the subgroups $K_{\mu^i, \sigma^i(\mu^j)}$ as in the Lemma \ref{lem:DefiAdmissible}. As a finite intersection of bounded open subgroups, $K_{\un{\mu}}$ has these properties, too. For any $b \in H_{\un{\mu}, k}(G)$ and $g \in K_{\un{\mu}}$ we have as $g \in G(O_L)$
\[\mu(bg\sigma) = \mu(b\sigma) = \mu^1.\]
Thus assume by induction $\mu((bg\sigma)^j) = \mu^j$ for any $j < i$. Then
\begin{align*}
 \mu((bg\sigma)^i) & = \mu(b \cdot g \cdot \sigma((bg\sigma)^{i-1})) \qquad {\rm with \;} g \in K_{\mu^1, \sigma(\mu^{i - 1})} \\
  & = \mu(b \cdot \sigma((bg\sigma)^{i-1})) \\
  & = \mu((b\sigma(b)) \cdot \sigma(g) \cdot \sigma^2((bg\sigma)^{i-2}) \qquad {\rm with \;} \sigma(g) \in K_{\mu^2, \sigma^2(\mu^{i - 2})} \\
  & = \mu((b\sigma)^2 \cdot (bg\sigma)^{i-2})) \\
  & = \ldots \\
  & = \mu((b\sigma)^i) = \mu^i
\end{align*}
and indeed $b\cdot g \in H_{\un{\mu}, k}(G)$.
\exit

\rem{}{
In general there is no constant $N$ (not even depending on the first Hodge point) such that $\mu((b\sigma)^{N+1})$ is determined by $\mu((b\sigma)^i)$ for $i = 1, \ldots, N$. In particular the partition into the sets $H_{\un{\mu}, k}(G)$ will not get stationary for growing $k$. This phenomena already appears for $G = GL_2$, although the computations in section \ref{subsec:ExplicitGL2} show that (for $k \geq 2$) this can happen only in the case where $H_{\un{\mu}, k}(GL_2)$ is contained in a Newton stratum with basic Newton point.
}

\subsection{Comparison between Hodge resp. Newton points under field extensions}\label{subsec:BaseChange}
At several points especially in section \ref{sec:Explicit} we will need to pass to field extensions of $L$. We will now explain in which cases we may do so. Although no additional conditions are necessary for $G = GL_n$, we have to impose further restriction concerning Hodge points in the general setting as special points do not behave very well under field extensions. \\
In this section the index at the map giving Hodge or Newton points indicates the field we use to define the points in the reductive group.

\prop{\label{prop:BaseChangeGL}}{}{
Let $G = GL_n$ and $T \subset GL_n$ a maximal split torus. Let $L \subset L'$ be a finite field extension with ramification index $e \in \M{N}$. Finally let $b \in GL_n(L)$ some matrix. Then $b$ defines canonically an element $b' \in GL_n(L')$, i.e. $b'$ is the matrix with the same entries as $b$, which are now viewed as elements in $L'$. Then
\[\nu_L(b) = \frac 1e \nu_{L'}(b')\]
\[\mu_L(b) = \frac 1e \mu_{L'}(b')\]
where each value in the $n$-tuple of $\nu_{L'}(b')$ resp. $\mu_{L'}(b')$ is divided by $e$.
}

\prooof
Wlog. let $T$ be the diagonal torus. For the Newton point choose $c \in GL_n(L)$ such that $cb\sigma(c)^{-1}$ is the normal form of $b$ (cf. Remark \ref{rem:GLNormalForm}) over $L$. Identifying both matrices $b$ and $c$ with elements in $GL_n(L')$ the same equation computes the Newton slope of $b' \in GL_n(L')$. The only difference is that the matrix coefficients $\unif^h$ in the normal form are no longer the $h$-th power but the $eh$-th power of a uniformizer of the ground field. Hence the Newton slopes of $b' \in GL_n(L')$ are the $e$th multiples of the slopes of $b \in GL_n(L)$. \\
For Hodge points choose $c_1, c_2 \in GL_n(O_L)$ such that $c_1bc_2 = {\rm diag}(\unif^{\mu(b)_1}, \ldots, \unif^{\mu(b)_n})$. Again the same computation gives the Hodge point of $b' \in GL_n(L')$ but with the difference as above that $\unif \in L$ is no longer a uniformizer of $L'$. 
\exit

\prop{\label{prop:BaseChangeGeneral}}{}{
Consider again the general setting as defined in \S \ref{subsec:Setting}. Let $L \subset L'$ be a finite field extension with ramification index $e \in \M{N}$  and fix an element $b \in G(L)$. This $b$ defines canonically an element $b' \in G(L')$. Then
\[\nu_L(b) = \frac 1e \nu_{L'}(b') \quad \in X_*(T)_{\M{Q}} / \Omega.\]
If there is a maximal $L$-split torus $S_L \subset T_L$ and a maximal $L'$-split torus $S_{L'} \subset T_{L'}$ already defined over $L$ and containing $S_L$, then identify the apartment $A(G_L, S_L, L)$ with the invariants of the apartment $A(G_{L'}, S_{L'}, L')$ under the Galois group (cf. \cite[\S1.10]{Tits}). Assume furthermore that the special point $\Lambda_0$ over $L$ stays special when considered in $A(G_{L'}, S_{L'}, L')$. Then
\[\mu_L(b) = \frac 1e \mu_{L'}(b') \quad \in X_*(T)_{\M{Q}} / \Omega.\]
}

\prooof
Consider first the case of Newton points. Let $\M{D}$ be the pro-algebraic torus over $\M{Q}_p$ with character group $\M{Q}$. Then we may consider $\nu_L(b) \in {\rm Int}\; G(L) \backslash {\rm Hom}_L(\M{D}, G)$ and $\nu_{L'}(b') \in {\rm Int} \; G(L') \backslash {\rm Hom}_{L'}(\M{D}, G)$ as the identification of ${\rm Int} \; G(L) \backslash {\rm Hom}_L(\M{D}, G)$ with $X_*(T)_{\M{Q}} / \Omega$ is functorial. \\
Now choose $\nu \in {\rm Hom}_L(\M{D}, G)$, $n > 0$ and $c \in G(L)$ s.th. the following three conditions hold for some uniformizer $\unif \in L$: \\
$i)$ $n\nu \in {\rm Hom}_L(\M{G}_m, G)$ \\
$ii)$ ${\rm Int}(c) \circ (n\nu)$ is defined over the fixed field of $\sigma^{n}$ on $L$ \\
$iii)$ $c (b\sigma)^n c^{-1} = c \cdot (n\nu)(\unif) \cdot c^{-1} \cdot \sigma^n$ \\
Note that $\nu$ and $\nu_L(b)$ coincide by definition in the quotient ${\rm Int} \; G(L) \backslash {\rm Hom}_L(\M{D}, G)$. \\
We will now check that $\nu' = e \cdot \nu$ (where we view $\nu$ via base-change as an element in ${\rm Hom}_{L'}(\M{D}, G)$) has the same properties with respect to $L'$: For this let $n' = n$ and $c'$ the $L'$-valued point defined by $c$. Then \\
$i)$ As $n\nu$ lies in the integral cocharacter group of $G$, its base-change to $L'$ lies in the integral cocharacter group of $G'$. But $n' \nu' = n \cdot (e\nu) = e \cdot (n\nu)$. Hence $n'\nu' \in {\rm Hom}_{L'}(\M{G}_m, G)$. \\
$ii)$ ${\rm Int}(c') \circ (n'\nu')$ is the $e$-th power of the base-change of ${\rm Int}(c) \circ (n\nu)$ to $L'$. Thus it is defined over the fixed field of $\sigma^n$ on $L'$. \\
$iii)$. First note that for some uniformizer $\unif' \in L'$ we have:
\[(n'\nu')(\unif') = (n\nu)(\unif'^e) = (n\nu)(\unif)\]
(the rightmost element considered via $G(L) \subset G(L')$). Thus we have
\[c' (b'\sigma)^{n'} c'^{-1} = c' \cdot (n'\nu')(\unif') \cdot c'^{-1} \cdot \sigma^{n'}\]
as both sides lie in $G(L)$ and are equal there. \\
This shows that $\nu'$ is exactly $\nu_{L'}(b')$ as desired. \vspace{1mm} \\
The case of Hodge points is easier: As $val_{L'}|_L = e \cdot val_L$ for the valuations, we get by definition 
\[v_L = \frac 1e v_{L'}|_{Z_L(L')}: Z_L(L') \to X_*(T)_{\M{Q}}\]
(denoting the centralizer of $S_L$ by $Z_L$). As the Hodge points are essentially the map $v_L$ respectively $v_{L'}$ one has:
\[\mu_L(b) = \frac 1e \mu_{L'}(b')\]
\exit

\rem{}{
i) Note that we change only the field $L$ and do not pass to some base-changed group. \\
ii) The additional assumptions for Hodge points are satisfied for unramified groups $G$ with hyperspecial vertices: As $G$ is quasi-split, the assumption on the existence of $S_{L'}$ is automatic. Furthermore hyperspecial vertices are by definition stable under finite field extensions.
}

\section{Interplay of \texorpdfstring{$\mathbf{H_{\un{\mu}, N}(G)}$}{refined Hodge strata} and Newton points for \texorpdfstring{$\mathbf{G = GL_n}$}{GL(n)}} \label{sec:GLn}
The main result of this section, Theorem \ref{thm:GLMain}, states that for sufficiently large $N$ (depending on $GL_n$ and the first Hodge point in $\un{\mu}$) the stratum $H_{\un{\mu}, N}(GL_n)$ lies inside one Newton stratum, i.e. the fixed Hodge points determine the Newton point. \\
To show this we will find an element $b_{min, \nu} \in GL_n(L)$ for each $b \in H_{\un{\mu}, N}(GL_n)$ such that one can bound the two differences 

i) $|\nu((b_{\min, \nu}\sigma)^k), \mu((b_{\min, \nu}\sigma)^k)|$ 

ii) $|\mu((b_{\min, \nu}\sigma)^k), \mu((b\sigma)^k)|$ \\
independently of $k$. \vspace{5mm} \\
Throughout this paragraph fix $|\cdot ,  \cdot|$ as in Remark \ref{rem:MetricGLn}, i.e. for elements with the same endpoint of their polygons:
\[|(\lambda_1, \ldots, \lambda_n), (\lambda_1', \ldots, \lambda_n')| = \sum_{\lambda_i > \lambda_i'} (\lambda_i - \lambda_i') \]
where $\sum_{i = 1}^n \lambda_i = \sum_{i = 1}^n \lambda_i'$ and $\lambda_1 \geq \ldots \geq \lambda_n$, $\lambda_1' \geq \ldots \geq \lambda_n'$ are in dominant order.

\subsection{Comparison for minimal elements}\label{sec:GLnMinimalElements}
In this section we will bound $|\nu((b_{\min, \nu}\sigma)^k), \mu((b_{\min, \nu}\sigma)^k)|$. The $b_{\min, \nu}$ will be the minimal elements in a Newton stratum corresponding to the minimal $p$-divisible groups as introduced in \cite{OortMinimal}. 

\defi{}{}{
i) A Newton point $\nu = (\frac hm, \ldots, \frac hm) \in \M{Q}^m$ is called superbasic if $gcd(h, m) = 1$. \\
ii) Let $\nu$ be superbasic with slope $\frac hm$. Let $e_0, \ldots, e_{m-1}$ be the standard basis of $L^m$ and define inductively $e_{k + m} = \unif e_k$ (for $k \in \M{Z}$). Then define $b_{min, \nu} \in GL_m(L)$ by $b_{min, \nu}(e_k) = e_{k + h}$ for $k = 0, \ldots, m - 1$. \\
iii) Let $\nu = (\nu_1, ..., \nu_n)$ be an arbitrary Newton point with dominantly ordered slopes. Let $\nu^1 = (\nu_1, \ldots, \nu_{i_1})$, $\nu^2 = (\nu_{i_1 + 1}, \ldots, \nu_{i_2})$, $\ldots$, $\nu^j = (\nu_{i_{j-1} + 1}, \ldots, \nu_n)$ be its superbasic parts. Then define $b_{min, \nu}$ to be the block matrix
\[b_{min, \nu} = \lmat{cccc} b_{min, \nu^1} & 0 & \ldots & 0 \\ 0 & b_{min, \nu^2} & \ldots & 0 \\ \vdots & \vdots & \ddots & \vdots \\ 0 & 0 & \ldots & b_{min, \nu^j} \rmat\]
The element $b_{min, \nu}$ is called a minimal element in the Newton stratum corresponding to $\nu$.
}

\rem{}{
In general a $\sigma$-conjugacy class $[b] \in B(G)$ is called superbasic if no representative of $[b]$ lies in a proper Levi subgroup defined over $F$ (cf. \cite[\S5.9]{GHKR}). In the case $G = GL_n$ this happens if the Newton point of $[b]$ is superbasic according to the definition given above. 
}

\lem{}{}{
i) $b_{min, \nu}$ has indeed Newton point $\nu$. \\
ii) $(b_{min, \nu}\sigma)^k = b_{min, k\nu}$ if $gcd(k, m) = 1$. Otherwise one can conjugate $(b_{min, \nu}\sigma)^k$ by a $\sigma$-invariant permutation matrix into $b_{min, k\nu}$.
}

\prooof
Obviously it suffices to check both statements for superbasic Newton points $\nu$. \\
i) As $gcd(h, m) = 1$ consider the basis $e'_i = e_{hi}$ for $i = 0, \ldots, m - 1$. Then $b_{min, \nu}(e'_i) = e'_{i + 1}$ for $i = 0, \ldots, m-2$ and $b_{min, \nu}(e'_{m - 1}) = e_{hm} = \unif^h e'_0$. The base-change matrix from $(e_i)$ to $(e'_i)$ is $\sigma$-invariant as it has entries $0$ or $\unif^k$ ($k \in \M{Z}$). Thus one can $\sigma$-conjugate $b_{min, \nu}$ into the standard form for the basic Newton-point $\nu$. \\
ii) $b_{min, \nu}$ is $\sigma$-invariant, thus 
\[(b_{min, \nu}\sigma)^k(e_i) = b_{min, \nu}^k(e_i) = e_{i + kh}.\]
If $gcd(k, m) = 1$, this is the definition of $b_{min, k\nu}$. \\
If $gcd(k, m) = s > 1$ and $t = \frac ms$, then let $\nu' = (\frac {kh}m, \ldots, \frac {kh}m) \in \M{Q}^t$ be superbasic with slope $\frac {kh}m$. Denote by $\{e'_j\}_{j \in \M{Z}}$ the basis vectors used in the definition of $b_{min, \nu'}$. For fixed $i = 0, \ldots s - 1$ we have a well-defined morphism $\varphi_i: L^t \to L^m$ of vector spaces with $\varphi_i(e'_j) = e_{i + sj}$. Then 
\[b_{min, \nu}^k(\varphi_i(e'_j)) = b_{min, \nu}^k(e_{i + sj}) = e_{i + sj + kh} = \varphi_i(e'_{j + \frac {kh}s}) = \varphi_i(b_{min, \nu'}(e'_j)).\]
Thus the restriction of $b_{min, \nu}$ to each subspace $\langle e_{i + sj} \rangle_{j = 0, \ldots, t-1}$ is via $\varphi_i$ equal to $b_{min, \nu'}$. Hence after a permutation of the basis vectors $b_{min, \nu}^k$ and $b_{min, k\nu}$ are equal.
\exit

\prop{}{}{
For each $b_{min, \nu} \in GL_n(L)$ and $k \in \M{N}$ we have:
\[|\nu((b_{min, \nu}\sigma)^k), \mu((b_{min, \nu}\sigma)^k)| \leq \frac n4\]
}

\prooof
By part ii) of the previous lemma, it suffices to show this for $k = 1$. Note for this that a permutation of the basis vectors does not change $|\cdot,\cdot|$. \\
Let us first consider the superbasic case. Then the Hodge slopes are exactly $\left\lfloor \frac {i + h}n\right\rfloor$ for $i = n-1, n-2, \ldots, 0$. It follows that if $\frac hn = a + \frac bn$ with $a \in \M{Z}$, $0 \leq b < n$, $b_{min, \nu}$ has Hodge slopes \[\mu_i(b_{min, \nu}) = \left\{\begin{array}{ccc} a + 1 & {\rm if} & i \leq b \\ a  & {\rm if} & i > b \end{array}\right.\]
Thus to compute $|\nu(b_{min, \nu}), \mu(b_{min, \nu})|$ one has to sum up the differences of the first $b$ slopes:
\[|\nu(b_{min, \nu}), \mu(b_{min, \nu})| = \left|b(a + 1) - b \frac hn\right| = \left|b - \frac {b^2}n\right| = \left|\frac n4 - \frac 1n(\frac n2 - b)^2\right| \leq \frac n4\]
Let now $\nu = (\nu_1, \ldots, \nu_n)$ be the dominant representative of an arbitrary Newton point. Then the slopes of $\mu(b_{min, \nu}) = (\mu_1, \ldots, \mu_n)$ can be permuted by some $\omega \in \Omega = S_n$ such that for each superbasic part $(\nu_i, \nu_{i + 1}, \ldots, \nu_{j})$ the tuple $(\mu_{\omega i}, \ldots, \mu_{\omega j})$ is the Hodge point of the corresponding superbasic minimal element (again dominantly ordered). Let now $I$ be the index set where $\mu_i \geq \nu_i$ and $J$ be the index set where $\mu_{\omega i} \geq \nu_i$. First note that $\omega$ preserves the sets $\{i \in \{1, 2, \ldots, n\} | k \leq \nu_i < k + 1\}$. Thus $\mu_{\omega^{-1} j} - \nu_{\omega^{-1} j} \geq 0$ if $\mu_{\omega^{-1} j} - \nu_j \geq 0$. As the converse is obvious we have $I = \omega J$ and see using $\nu(b_{min, \nu}) \prec \mu(b_{min, \nu})$:
\[|\nu(b_{min, \nu}), \mu(b_{min, \nu})| = \sum_{i \in I} (\mu_i - \nu_i) \leq \sum_{i \in I} (\mu_i - \nu_{\omega^{-1} i}) = \sum_{j \in J} (\mu_{\omega j} - \nu_j) \leq \frac n4\]
where the last inequality follows from the estimate for superbasic Newton points applied to each superbasic part of $\nu$.
\exit

\subsection{Bounding the \texorpdfstring{$\mathbf{\sigma}$}{Frobenius}-conjugating elements}\label{sec:GLnSigmaConjugation}
Now we will show that any two elements in $GL_n$ with the same Newton and Hodge point are $\sigma$-conjugates by some element with bounded Hodge point. To do so we use the main result of \cite{RapoZinkBTbuilding} in the case of $G = GL_n$. In this case it can be seen as a direct consequence of the Rapoport-Zink lemma (cf. \cite[2.17-2.19]{RZ1} or \cite[prop. 1.6]{RapoZinkBTbuilding}). The general theorem will be stated in \S \ref{subsec:GeneralFinal}. 

\defi{}{}{
Let $G$ be any reductive group defined over $F$, $b \in G(L)$. Then let $J_b$ be the connected linear algebraic group over $F$ with 
\[J_b(R) = \{g \in G(R \otimes_F L) \; | \; \sigma(g) = b^{-1}gb\}\]
for any $F$-algebra $R$.
}

\rem{}{
That $J_b$ is indeed representable as a connected linear algebraic group is shown in \cite[prop. 1.12]{RZ1}.
}

\defi{\label{defi:GLDecencyEquation}}{(cf. \cite[def. 1.8]{RZ1} or loc. cit. def. \ref{defi:DecencyEquation})}{
An element $b \in GL_n(L)$ is said to fulfill a \textit{decency equation} if there is an $s \in \M{N}$ and a uniformizer $\unif \in L$ with
\[(b\sigma)^s = {\rm diag}(\unif^{s\nu_1}, \unif^{s\nu_2}, \ldots, \unif^{s\nu_n})\] 
with Newton slopes $\nu_1, \ldots, \nu_n$ of $b$ (not necessarily in dominant order). 
}

\rem{}{
Minimal elements always fulfill decency equations.
}

\prop{}{(cf. \cite[thm. 1.4]{RapoZinkBTbuilding} for $G = GL_n$)}{
Let $b_0 \in GL_n(L)$ fulfilling a decency equation relative to $s \in \M{N}$ and fix $\mu \in X_*(T)_{\M{Q}} / \Omega$. Then there is a finite set $\MC{S} \subset X_*(T)_{\M{Q}} / \Omega$ with the following property: For any $g \in GL_n(L)$ with $\mu(g^{-1}b_0\sigma(g)) = \mu$ there is a $j \in J_{b_0}(F_s)$ with $\mu(jg) \in \MC{S}$. \\
Here $F_s$ denotes the unramified extension of degree $s$ of $F$. 
}

\prop{\label{prop:GLBoundingSigmaConjugacy}}{}{
Fix $\nu \in X_*(T)_{\M{Q}} / \Omega$, $\mu \in X_*(T)_{\M{Q}} / \Omega$. Then there is a finite set $\MC{S} \subset X_*(T)_{\M{Q}} / \Omega$ such that for any $b \in GL_n(L)$ with $\mu(b) = \mu$ and $\nu(b) = \nu$ there is a $g \in GL_n(L)$ and a minimal element $b_{min, \nu}$ with $\mu(g) \in  \MC{S}$ and $b = g^{-1} \cdot b_{min, \nu} \cdot \sigma(g)$.
}

\prooof
Choose a finite set $\MC{S} \subset X_*(T)_{\M{Q}} / \Omega$ for $b_0 = b_{min, \nu}$ as in the previous proposition. Let now $b \in GL_n(L)$ with $\mu(b) = \mu$ and $\nu(b) = \nu$. As all elements are conjugate to some minimal element there is some $g' \in GL_n(L)$ and a minimal element $b_{min, \nu}$ with $b = g'^{-1} \cdot b_{min, \nu} \cdot \sigma(g')$. Then $g'$ fulfills all conditions of the previous proposition and there is a $j \in J_{b_{min, \nu}}(F_s) \subset GL_n(L)$ with $g := jg' \in \MC{S}$. Finally note that by definition of $J_{b_{min, \nu}}$
\[g^{-1} \cdot b_{min, \nu} \cdot \sigma(g) = g'^{-1} \cdot j^{-1} \cdot b_{min, \nu} \cdot \sigma(j) \cdot \sigma(g') = g'^{-1} \cdot b_{min, \nu} \cdot \sigma(g') = b.\] 
\exit

\subsection{Final estimates}\label{subsec:GLnFinal}
The main remaining part is to show that for $b, g \in GL_n(L)$ one can give bounds on the Hodge points of $g^{-1}b\sigma(g)$ in terms of the Hodge point of $b$ and $g$. 

\defi{}{}{
Let $G$ be a reductive group. For $x, x' \in X_*(T)_{\M{R}} / \Omega$ let $\tilde{x}$ resp. $\tilde{x}'$ be the dominant representatives of $x$ resp. $x'$ in $\ov{C} \subset X_*(T)_{\M{R}}$. Then let \\
i) $x \oplus x'$ be the image of $\tilde{x} + \tilde{x}' \in X_*(T)_{\M{R}}$ in the quotient $X_*(T)_{\M{R}} / \Omega$. \\
ii) $x \oplus_{\omega_0} x'$ be the image of $\tilde{x} + \omega_0 \tilde{x}' \in X_*(T)_{\M{R}}$ in the quotient $X_*(T)_{\M{R}} / \Omega$, where $\omega_0 \in \Omega$ is the longest element.
}

\rem{}{
In the case of $G = GL_n$, $\tilde{x} = (\lambda_1, \ldots, \lambda_n)$, $\tilde{x}' = (\lambda_1', \ldots, \lambda_n')$ (both dominant)
\begin{align*}
 x \oplus x'            & = (\lambda_1 + \lambda_1', \lambda_2 + \lambda_2', \ldots, \lambda_n + \lambda_n') \\
 x \oplus_{\omega_0} x' & = (\lambda_1 + \lambda_n', \lambda_2 + \lambda_{n - 1}',  \ldots, \lambda_n + \lambda_1')
\end{align*}
Note that the last element need not lie in the dominant chamber.
}

\lem{\label{lem:GLHodgeEstimate}}{}{
Let $b, b' \in GL_n(L)$. Then
\[\mu(b) \oplus_{\omega_0} \mu(b') \prec \mu(bb') \prec \mu(b) \oplus \mu(b').\]
}

\prooof
We will first consider only the smallest Hodge slope: To simplify the situation use the Cartan decomposition to write $b' = c_1b_0'c_2$ with $b_0' = {\rm diag}(\unif^{\mu_1'}, \ldots, \unif^{\mu_n'}) \in T(L)$, $\mu_1' \geq \ldots \geq \mu_n'$, $c_1, c_2 \in GL_n(O_L)$. Now decompose $bc_1 = c_3b_0c_4$ with $c_3 \in GL_n(O_L)$, $b_0 = {\rm diag}(\unif^{\mu_{\tau(1)}}, \ldots, \unif^{\mu_{\tau(n)}}) \in T(L)$, $\tau \in S_n$ a permutation, $\mu_1 \geq \ldots \geq \mu_n$, $c_4 \in U(O_L)$ (where $U \subset B$ denotes the unipotent radical). Then replacing $b$ by $b_0c_4$ and $b'$ by $b_0'$ does not change any of the considered Hodge points. The algorithm to compute elementary divisors implies, that the smallest Hodge slope equals the smallest valuation of a matrix coefficient of $b_0c_4b_0'$. Hence (denoting the coefficients of $c_4$ by $(c_4)_{i,j}$)
\[\mu(bb')_n = \mu(b_0c_4b_0')_n = \min_{i, j \in \{1, \ldots, n\}} (\mu_{\tau(i)} + \mu_j' + v((c_4)_{i,j})).\]
In particular one has
\[\mu(bb')_n \geq \mu_n + \mu_n'.\]
For the upper bound note that for each $j \in \{1, \ldots, n\}$ there is a $i \in \{1, \ldots, n\}$ with $j \leq i$ and $n - j + 1 \leq \tau(i)$ (otherwise the $n - j + 1$ elements $j, \ldots, n$ would be mapped under $\tau$ to the $n - j$ elements $1, \ldots, n - j$). Hence
\[\mu_{n - j + 1} + \mu_j' \geq \mu_{\tau(i)} + \mu_i' \geq \mu(bb')_n\] 
and $\mu(bb')_n \leq \min_i (\mu(b) \oplus_{\omega_0} \mu(b'))_i$. \vspace{2mm} \\
For further Hodge slopes let $GL_n \cong GL(V)$ for some vector space $V$ and consider $b^{\wedge k} \in GL(\bigwedge^kV)$. If the slopes of $b$ are $\mu_1 \geq \ldots \geq \mu_n$ (as above) then the slopes of $b^{\wedge k}$ are $\{\sum_{i \in S} \mu_i \; | \; S \subseteq \{1, 2, \ldots, n\}, |S| = k\}$. Similarly for $b'$ and $bb'$. Thus we get for all $k = 1, \ldots, n$
\begin{align*}
\sum_{i = n-k+1}^n \mu(bb')_i & = \mu((bb')^{\wedge k})_{\binom{n}{k}} = \mu(b^{\wedge k}b'^{\wedge k})_{\binom{n}{k}} \geq \mu(b^{\wedge k})_n + \mu(b'^{\wedge k})_n  \\
& = \sum_{i = n-k+1}^n \mu(b)_i + \sum_{i = n-k+1}^n \mu(b')_i = \sum_{i = n-k+1}^n (\mu(b) \oplus \mu(b'))_i
\end{align*}
thus the right hand side of the inequality. \\
Similarly fix any $S \subset \{1, \ldots n\}$, $|S| = k$. Then by an analogue argument as for the smallest slope, there is a $j \in \{1, \ldots, \binom{n}{k}\}$ such that
\[\sum_{i \in S} \mu(b)_i \geq \mu(b^{\wedge k})_j \qquad {\rm and} \qquad \sum_{i \in S} \mu(b')_{n - i + 1} \geq \mu(b'^{\wedge k})_{\binom{n}{k} - j + 1}.\]
Hence
\begin{align*}
\sum_{i \in S} \mu(b)_i + \sum_{i \in S} \mu(b')_{n - i + 1} & \geq  \mu(b^{\wedge k})_j + \mu(b'^{\wedge k})_{\binom{n}{k} - j + 1}  \\
& \geq  \min_{i \in \left\{1, \ldots, \binom{n}{k} \right\}} (\mu(b^{\wedge k})_i + \mu(b'^{\wedge k})_{\binom{n}{k} - i + 1}) \\
& \geq  \mu(b^{\wedge k}b'^{\wedge k})_{\binom{n}{k}} = \mu((bb')^{\wedge k})_{\binom{n}{k}} = \sum_{i = n-k+1}^n \mu(bb')_i.
\end{align*}
That means $\sum_{i = n-k+1}^n \mu(bb')_i$ is smaller that the sum of the $k$ least elements in $\mu(b) \oplus_{\omega_0} \mu(b')$. This gives the left hand side of the inequality.
\exit

\rem{}{
When only considering $G = GL_n$ the right hand side of the inequality would suffice for our purposes. Only in the next chapter we will really need both sides to derive similar statements for arbitrary connected reductive groups. 
}

\prop{\label{prop:GLHodgeForMinimal}}{}{
Let $\nu \in X_*(T)_{\M{Q}} / \Omega$ be any Newton point and $\mu \in X_*(T)_{\M{Q}} / \Omega$ any Hodge point. Then there is a constant $C$ such that for any $b \in GL_n(L)$ with $\nu(b) = \nu$ and $\mu(b) = \mu$ and any $k \in \M{N}$
\[|\mu((b\sigma)^k), \mu((b_{min, \nu}\sigma)^k)| < C.\]
}

\prooof
$b_{min, \nu}$ satisfies a decency equation and we may apply Proposition \ref{prop:GLBoundingSigmaConjugacy}. Thus we can write $b = g^{-1} \cdot b_{min, \nu} \cdot \sigma(g)$ for some $g \in GL_n(L)$ with Hodge point in some finite set $\MC{S}$. Applying Proposition \ref{lem:GLHodgeEstimate} gives together with $\mu(g) = \mu(\sigma(g))$
\begin{align*}
 & \mu(g^{-1}) \oplus_{\omega_0} \mu((b_{min, \nu}\sigma)^k) \oplus_{\omega_0} \mu(g) \\
 & \qquad = \mu(g^{-1}) \oplus_{\omega_0} \mu((b_{min, \nu}\sigma)^k) \oplus_{\omega_0} \mu(\sigma^k(g)) \\
 & \qquad \prec \mu(g^{-1} \cdot (b_{min, \nu}\sigma)^k \cdot g) = \mu((g^{-1}b_{min, \nu}\sigma(g)\sigma)^k) =  \mu((b\sigma)^k) \\
 & \qquad \prec \mu(g^{-1}) \oplus \mu((b_{min, \nu}\sigma)^k) \oplus \mu(\sigma^k(g))  \\
 & \qquad = \mu(g^{-1}) \oplus \mu((b_{min, \nu}\sigma)^k) \oplus \mu(g). 
\end{align*}
$\mu(g)$ and $\mu(g^{-1})$ are bounded because $\MC{S}$ is finite. Hence there is some constant $C > 0$ with
\[|\mu((b\sigma)^k), \mu((b_{min, \nu}\sigma)^k)| < C.\]
\exit

\prop{\label{prop:GLDifferenceNewtonPolygon}}{}{
For any two $b, b' \in GL_n(L)$ with $\nu(b) \neq \nu(b')$, there is some value at which the (concave) Newton polygons of $\nu(b)$ and $\nu(b')$ differ by at least $\frac 1n$. In particular
\[|\nu(b), \nu(b')| \geq \frac 1n.\]
}

\prooof
As $\nu(b) \neq \nu(b')$ one can find at least one vertex of the Newton polygon of either $\nu(b)$ or $\nu(b')$ which does not lie on the Newton polygon of the other. At this point, the difference of both polygons is at least $\frac 1n$. But $|\nu(b), \nu(b')|$ is at least as big as any vertical difference between those polygons.
\exit

\thm{\label{thm:GLMain}}{}{
Let $\mu^1$ be a Hodge point. Then there is a constant $C'$ only depending on $\mu^1$, such that each stratum $H_{\un{\mu}, N}$ with $N > C'$ and $\un{\mu} = (\mu^1, \mu^2, \ldots, \mu^N)$ lies inside a Newton stratum. 
}

\prooof
Let $\mu$ be a Hodge point. Fix a Newton point $\nu$ for now. Let $b \in GL_n(L)$ with Newton point $\nu$ and Hodge point $\mu$. Then with the constant $C$ of Proposition \ref{prop:GLHodgeForMinimal}
\begin{align*}
|k\nu(b), \mu((b\sigma)^k)| & =  |\nu((b\sigma)^k), \mu((b\sigma)^k)| \\
& \leq  |\nu((b\sigma)^k), \mu((b_{min, \nu}\sigma)^k)| + |\mu((b_{min, \nu}\sigma)^k), \mu((b\sigma)^k)| \\
& \leq  \frac n4 + C 
\end{align*}
and for $N > C'_{\nu} := n(\frac n4 + C)$:
\[\left|\nu(b), \frac 1N \mu(b\sigma)^N\right| \leq \frac 1N (\frac n4 + C) < \frac 1n.\] 
Thus the (concave) Newton polygon of $\nu(b)$ lies below the polygon of $\frac 1N \mu(b\sigma)^N$ by Mazur's inequality and differs from it by less than $\frac 1n$ at every point. Hence by Proposition \ref{prop:GLDifferenceNewtonPolygon} $\mu(b\sigma)^{C'_{\nu}}$ determines $\nu(b)$.
Now let $C' = \max_{\nu \prec \mu^1} (C'_{\nu})$ and fix some $\un{\mu} = (\mu^1, \ldots, \mu^N)$ for $N > C'$. Let $b, b' \in H_{\un{\mu}, N}$. Then $\nu(b), \nu(b') \prec \mu^1$ by Mazur's inequality. By definition of $C'$, the Hodge point $\mu^N$ determines both $\nu(b)$ and $\nu(b')$. Thus $\nu(b) = \nu(b')$, i.e. $H_{\un{\mu}, N}$ lies inside a Newton stratum.
\exit

\section{Interplay of \texorpdfstring{$\mathbf{H_{\un{\mu}, n}}$}{refined Hodge strata} and Newton points for general \texorpdfstring{$\mathbf{G}$}{G}} \label{sec:GeneralG}
We use the same strategy as for $GL_n$ in the general context. Nevertheless additional arguments are needed at several points.

\subsection{Comparison for chosen elements}\label{subsec:GeneralStandardElements}
There is no well-established notion of 'minimal elements' for arbitrary connected reductive groups (although there is some generalization under additional assumptions). Hence we use some generalization of the normal form of elements in $GL_n$ instead.

\lem{\label{lem:GeneralHodgeBound}}{}{
Let $b, b' \in G(L)$. Then
\[\mu(b) \oplus_{\omega_0} \mu(b') \prec \mu(bb') \prec \mu(b) \oplus \mu(b').\]
}

\prooof
It suffices to check this for the images under all representations $G \to GL(V)$ (cf. \cite[lemma 2.2]{RapoRich}). But for $GL(V)$ this was shown in Lemma \ref{lem:GLHodgeEstimate}.
\exit

\prop{\label{prop:GeneralStandardEstimate}}{}{
For each $b \in G(L)$ there is a $b_0 \in G(L)$ in the $\sigma$-conjugacy class of $b$ and a constant $C_0 > 0$ such that for all $k \in \M{N}$ we have
\[|\nu(b_0\sigma)^k, \mu(b_0\sigma)^k| \leq C_0.\] 
}

\prooof
To any $b \in G(L)$ there is a representative $\nu \in {\rm Hom}_L(\M{D}, G)$ of the Newton point $\nu(b)$ (viewed as an element in ${\rm Int} \; G(L) \backslash {\rm Hom}_L(\M{D}, G)$), $n > 0$, $\unif \in L$ a uniformizer and an element $c \in G(L)$ such that
\[c (b\sigma)^n c^{-1} = c \cdot (n\nu)(\unif) \cdot c^{-1} \cdot \sigma^n.\]
Let $b_0 = cb\sigma(c)^{-1}$ and $\nu_0 = c \nu c^{-1}$. Then
\[(b_0\sigma)^n = (cb\sigma(c)^{-1}\sigma)^n = c (b\sigma)^n c^{-1} = c \cdot (n\nu)(\unif) \cdot c^{-1} \cdot \sigma^n = (n \nu_0)(\unif) \cdot \sigma^n.\]
Hence $(b_0\sigma)^n$ lies in some torus and trivially
\[\mu((b_0\sigma)^n) = \nu((b_0\sigma)^n) = \nu = \nu_0 \in {\rm Int} \; G(L) \backslash {\rm Hom}_L(\M{D}, G).\]
For general $k = xn + y$, $x \in \M{N}$, $0 \leq y < n$ we have
\begin{align*}
\mu((b_0\sigma)^{xn}) \oplus_{\omega_0} \mu((b_0\sigma)^y) & =  \mu((b_0\sigma)^{xn}) \oplus_{\omega_0} \mu(\sigma^{xn}(b_0\sigma)^y) \\
& \prec  \mu((b_0\sigma)^{xn + y})  \\
& \prec  \mu((b_0\sigma)^{xn}) \oplus \mu(\sigma^{xn}(b_0\sigma)^y)  \\
& =  \mu((b_0\sigma)^{xn}) \oplus \mu((b_0\sigma)^y). 
\end{align*}
Thus there is a constant $C_1 > 0$ such that for all $k = xn + y$
\[|\mu((b_0\sigma)^{xn + y}), \mu((b_0\sigma)^{xn}| < C_1.\] 
Furthermore as $\nu((b\sigma)^k) = k \nu(b\sigma)$ one has a constant $C_2$ with
\[|\nu((b_0\sigma)^{xn + y}), \nu((b_0\sigma)^{xn}| = |\nu((b_0\sigma)^y), 0| < C_2\] 
which finally allows us to conclude
\[|\nu((b_0\sigma)^{xn + y}), \mu((b_0\sigma)^{xn + y})| < C_1 + C_2 + |\nu((b_0\sigma)^{xn}), \mu((b_0\sigma)^{xn})| = C_1 + C_2 \eqqcolon C_0.\]
\exit

\rem{}{
Contrary to the case $G = GL_n$ this proof yields a constant depending on the $\sigma$-conjugacy class and not only on $G$ itself. 
}

\subsection{Bounding the \texorpdfstring{$\mathbf{\sigma}$}{Frobenius}-conjugating elements for general \texorpdfstring{$\mathbf{G}$}{G} and final estimates}\label{subsec:GeneralFinal}
Recall the definition of $J_b$, $F_s$ as given in section \ref{sec:GLnSigmaConjugation} and generalize Definition \ref{defi:GLDecencyEquation} as follows:

\defi{\label{defi:DecencyEquation}}{(\cite[def. 1.8]{RZ1})}{
An element $b \in G(L)$ is said to fulfill a \textit{decency equation} if there is a $s \in \M{N}$ and a uniformizer $\unif \in L$ with
\[(b\sigma)^s = s\nu_b(\unif) \cdot \sigma^s \; \in G(L) \rtimes \langle \sigma\rangle.\]
}
The main tool to bound the difference between Hodge points of $\sigma$-conjugated elements is again the main theorem of \cite{RapoZinkBTbuilding}, now in full generality:

\thm{}{(\cite[thm. 1.4]{RapoZinkBTbuilding})}{
Fix $b \in G(L)$ and $s > 0$ such that a decency equation holds for $b$ relative to $s$. For $c > 0$ there exists a bound $C > 0$ with the following property: If $x \in \MC{B}(G, L)$ is an element in the extended Bruhat-Tits building over $L$ such that $d(x, b\sigma(x)) < c$, then there exists a $x_0 \in \MC{B}(J, F_s)$ with $d(x, x_0) < C$. 
}

\prop{\label{prop:GeneralHodgeForStandard}}{}{
Let $\nu \in X_*(T)_{\M{Q}} / \Omega$ be any Newton point and $\mu \in X_*(T)_{\M{Q}} / \Omega$ any Hodge point. Then there is a constant $C$ such that for any $b \in G(L)$ with $\nu(b) = \nu$ and $\mu(b) = \mu$ and any $k \in \M{N}$
\[|\mu(b\sigma)^k, \mu(b_0\sigma)^k| < C.\]
}

\prooof
Copy the proof of Proposition \ref{prop:GLBoundingSigmaConjugacy} and \ref{prop:GLHodgeForMinimal}, replace $GL_n$ by $G$ and note that by definition $b_0$ satisfies a decency equation.
\exit

\thm{\label{thm:GeneralBound}}{}{
Let $\mu^1$ be a Hodge point of $G$. Then there is a constant $C'$ only depending on $\mu^1$, such that each stratum $H_{\un{\mu}, N}$ with $N > C'$ and $\un{\mu} = (\mu^1, \mu^2, \ldots, \mu^N)$ lies inside a Newton stratum. 
}

\prooof
Fix a Hodge point $\mu$. For now fix a Newton point $\nu$, too. Let $b \in G(L)$ with Newton point $\nu$ and Hodge point $\mu$. Then for $b_0 \in G(L)$ and $C_0$ as in Proposition \ref{prop:GeneralStandardEstimate} and $C$ as in Proposition \ref{prop:GeneralHodgeForStandard}
\begin{align*}
|k\nu(b), \mu((b\sigma)^k)| & =  |\nu((b\sigma)^k), \mu((b\sigma)^k)| \\
& \leq  |\nu((b\sigma)^k), \mu((b_0\sigma)^k)| + |\mu((b_0\sigma)^k), \mu((b\sigma)^k)| \leq C_0 + C.
\end{align*}
As there are only finitely many Newton points below $\mu$ (cf. \cite[prop. 2.4 iii)]{RapoRich}) one can find some constant $\varepsilon > 0$ such that for any two Newton points $\nu', \nu'' \prec \mu$ one has $|\nu', \nu''| > \varepsilon$. Then for $N > C'_{\nu} = 2\varepsilon^{-1}(C_0 + C)$
\[\left|\nu(b), \frac 1N \mu(b\sigma)^N\right| \leq \frac 1N (C_0 + C) < \frac{\varepsilon}2.\] 
Thus $\nu(b)$ is uniquely determined by any $\mu(b\sigma)^N$. \\
Now let $C' = \max_{\nu \prec \mu^1} (C'_{\nu})$ (again we take the maximum over a finite set) and fix some $\un{\mu} = (\mu^1, \ldots, \mu^N)$ for $N > C'$. Then for any two elements in $H_{\un{\mu}, N}$, the Hodge point $\mu^N$ determines the Newton point for both of them. Hence $H_{\un{\mu}, N}$ lies inside a Newton stratum.
\exit

\conj{}{}{
Let $G$ be a reductive group. Then there is a constant $C'$ depending only on $G$, such that each stratum $H_{\un{\mu}, N}$ with $N > C'$ and any $N$-tuple of Hodge points $\un{\mu}$ lies inside a Newton stratum.
}

\rem{}{
The proof of this conjecture needs a different approach. Even for $GL_2$ and any $k \in \M{N}$ there are matrices $b \in GL_2(L)$ where $\mu((b\sigma)^k)$ differs greatly from $k\nu(b)$. \\
However we prove the conjecture for $G = GL_2$ in Proposition \ref{prop:ExplicitGL2} and even for any scalar restriction of $GL_2$ in Proposition \ref{prop:ExplicitBoundResGL2}.
}

\section{Results for some explicit groups}\label{sec:Explicit}
In the last paragraph it was shown that $H_{\un{\mu}, N}$ lies in a Newton stratum for sufficiently large $N$. But in most cases the constants derived in the general proof are far from being optimal. This problem will be treated for some small groups to establish some small bounds, preferably $2$.

\subsection{\texorpdfstring{$\mathbf{GL_2}$}{GL(2)}}\label{subsec:ExplicitGL2}
\prop{\label{prop:ExplicitGL2}}{}{
For any pair of Hodge points $\un{\mu} = (\mu^1, \mu^2)$ of $GL_2(L)$ the stratum $H_{\un{\mu}, 2}(GL_2)$ lies in some Newton stratum.
}

\prooof
Let $b \in GL_2(L)$. If $b$ has non-integral Newton slopes, let $L'$ be a totally ramified extension of $L$ of degree $2$ and view $b$ as an element $b' \in GL_2(L')$. The Newton point of $b'$ over $L'$ lies now in the integral cocharacter group and we can recover the Newton point $\nu_L(b)$ out of the Newton point $\nu_{L'}(b')$ (cf. Proposition \ref{prop:BaseChangeGL}). Hence we may replace $L$ by $L'$ and $b$ by $b'$ and can assume wlog. that $\nu(b)$ lies in the (integral) cocharacter group. \\
Then there are elements $c \in GL_2(L)$, $\nu_1, \nu_2 \in \M{Z}$ with $\nu_1 \leq \nu_2$ and some uniformizer $\unif \in L$ such that
\[b = c \cdot \lmat{cc} \unif^{\nu_1} & 0 \\ 0 & \unif^{\nu_2}\rmat \cdot \sigma(c)^{-1}.\]
Applying the algorithm to compute elementary divisors to $c$ we can write
\[c = \lmat{cc} 1 & 0 \\ \lambda' \unif^{d'} & 1 \rmat \cdot \lmat{cc} 0 & 1 \\ 1 & 0\rmat^{e'} \cdot \lmat{cc} a_1 & 0 \\ 0 & a_2\rmat \cdot \lmat{cc} \unif^{x} & 0 \\ 0 & \unif^{y}\rmat \cdot \lmat{cc} 0 & 1 \\ 1 & 0\rmat^{e} \cdot \lmat{cc} 1 & 0 \\ \lambda \unif^{d} & 1\rmat\]
for some $x, y \in \M{Z}$, $y \leq x$, $d, d' \geq 0$, $e, e' \in \{0, 1\}$, $\lambda, \lambda' \in O_L^\times \cup \{0\}$ and $a_1, a_2 \in O_L^\times$. As neither Newton nor Hodge points change when $\sigma$-conjugating with some element of $GL_2(O_L)$, we may assume $\lambda' = 0$, $e' = 0$ and $a_1 = a_2 = 1$. Then
\begin{align*}
b & = \lmat{cc} \unif^{x} & 0 \\ 0 & \unif^{y}\rmat \cdot \lmat{cc} 0 & 1 \\ 1 & 0\rmat^{e} \cdot \lmat{cc} 1 & 0 \\ \lambda \unif^{d} & 1\rmat \cdot \lmat{cc} \unif^{\nu_1} & 0 \\ 0 & \unif^{\nu_2}\rmat \cdot \\
& \hspace{23mm} \cdot \lmat{cc} 1 & 0 \\ - \sigma(\lambda) \unif^{d} & 1\rmat \cdot \lmat{cc} 0 & 1 \\ 1 & 0\rmat^{e} \cdot \lmat{cc} \unif^{-x} & 0 \\ 0 & \unif^{-y}\rmat.
\end{align*}
\textbf{Case 1: $e = 0$} \\
Computing $b$ we get with $\delta = d - x + y$:
\[b = \lmat{cc} \unif^{\nu_1} & 0 \\ \lambda \unif^{\nu_1 + \delta} - \sigma(\lambda) \unif^{\nu_2 + \delta} & \unif^{\nu_2}\rmat\]
\[(b\sigma)^2 = \lmat{cc} \unif^{2\nu_1} & 0 \\ \lambda \unif^{2\nu_1 + \delta} - \sigma^2(\lambda) \unif^{2\nu_2 + \delta} & \unif^{2\nu_2}\rmat\]
If $\delta \geq 0$ or $\lambda = 0$ then $\mu(b) = (\nu_1, \nu_2)$, $\mu((b\sigma)^2) = (2\nu_1, 2\nu_2)$. \\
If $\delta < 0$ and $\lambda \neq 0$ we have to distinguish two further cases: \\
\textbf{Case 1.1: $\delta < 0$, $\lambda \neq 0$ and $\nu_1 \neq \nu_2$} \\
Then $\mu(b) = (\nu_1 + \delta, \nu_2 - \delta)$, $\mu((b\sigma)^2) = (2\nu_1 + \delta, 2\nu_2 - \delta)$. \\
\textbf{Case 1.2: $\delta < 0$, $\lambda \neq 0$ and $\nu_1 = \nu_2$} \\
Now let $\delta'_1, \delta'_2 \geq 0$ be the unique integers with 
\[\lambda - \sigma(\lambda) \in (\unif^{\delta_1'}) \setminus (\unif^{\delta_1' + 1}) \quad {\rm and} \quad \lambda - \sigma^2(\lambda) \in (\unif^{\delta_2'}) \setminus (\unif^{\delta_2' + 1})\]
Obviously $\delta_1' \leq \delta_2'$. \\
If $\delta + \delta_1' \geq 0$, then $\mu(b) = (\nu_1, \nu_2)$, $\mu((b\sigma)^2) = (2\nu_1, 2\nu_2)$. \\
If $\delta + \delta_1' < 0$ and $\delta + \delta_2' \geq 0$, then $\mu(b) = (\nu_1 + (\delta + \delta_1'), \nu_2 - (\delta + \delta_1'))$, $\mu((b\sigma)^2) = (2\nu_1, 2\nu_2)$. \\
If $\delta + \delta_2' < 0$, then $\mu(b) = (\nu_1 + (\delta + \delta_1'), \nu_2 - (\delta + \delta_1'))$, $\mu((b\sigma)^2) = (2\nu_1 + (\delta + \delta_2'), 2\nu_2 - (\delta + \delta_2'))$. \\ 
\textbf{Case 2: $e = 1$} \\
Computing $b$ again with $\delta = d + x - y \geq 0$:
\[b = \lmat{cc} \unif^{\nu_2} & \lambda \unif^{\nu_1 + \delta} - \sigma(\lambda) \unif^{\nu_2 + \delta} \\ 0 & \unif^{\nu_1}\rmat\]
\[(b\sigma)^2 = \lmat{cc} \unif^{2\nu_2} & \lambda \unif^{2\nu_1 + \delta} - \sigma^2(\lambda) \unif^{2\nu_2 + \delta} \\ 0 & \unif^{2\nu_1}\rmat\]
Thus $\mu(b) = (\nu_1, \nu_2)$, $\mu((b\sigma)^2) = (2\nu_1, 2\nu_2)$. \\
Hence in all cases we can recover the Newton point out of $\mu(b)$ and $\mu((b\sigma)^2)$ via the following procedure: \\
If $\mu(b) = (\mu_{1,1}, \mu_{1,2})$ and $\mu((b\sigma)^2) = (\mu_{2,1}, \mu_{2,2})$ are dominant representatives of the Hodge points, then compute $2\mu_{1, 1} - \mu_{2, 1}$. This value is negative if and only if we are in case $1.2$ and $\delta + \delta_1' < 0$. But then the Newton point is basic and can be computed as $\nu(b) = \frac 12(\mu_{1,1} + \mu_{1,2}, \mu_{1,1} + \mu_{1,2})$. But if $2\mu_{1, 1} - \mu_{2, 1} \geq 0$ (i.e. in all remaining cases), the formula $\nu(b) = (\mu_{2,1} - \mu_{1,1}, \mu_{2,2} - \mu_{1,2})$ holds. 
\exit

\rem{}{
The same case-by-case analysis but for $(b\sigma)^{n_1}$ and $(b\sigma)^{n_2}$ (with $n_2 > n_1 > 0$) instead of $b$ and $(b\sigma)^2$ shows that the Newton point can be recovered from $\mu((b\sigma)^{n_1})$ and $\mu((b\sigma)^{n_2})$ whenever $n_1 \,|\, n_2$. The last divisibility condition is necessary to ensure a similar inequality between the $\delta_i'$ as above. 
}

\subsection{Restriction of scalars for totally ramified extensions}\label{subsec:ExplicitRamified}
In this section we compare scalar restrictions for totally ramified extensions of some connected reductive group $G$ to the group itself. \\
For this comparison fix a tower of fields $\M{Q}_p \subset F' \subset F \subset \ov{\M{Q}}_p$ with $F'$ and $F$ finite over $\M{Q}_p$ and $F/F'$ totally ramified. Let as in \S \ref{subsec:Setting} $L = F \cdot K$ and $L' = F' \cdot K$. Then $L/L'$ is again a totally ramified extension of the same degree as $F/F'$.  \\
Let $G$ be a connected reductive group over $F$. Then by definition of scalar restrictions there is an isomorphism of abstract groups
\[\phi: \left({\rm Res}_{F/F'}(G)\right)(L') \cong G(F \otimes_{F'} L') = G(L).\]
Fix as usual a maximal torus $T \subset G$. Then $T' = {\rm Res}_{F/F'}(T) \subset {\rm Res}_{F/F'}(G)$ is again a maximal torus. We will denote their base-change to $L$ resp. $L'$ by $T_L$ resp. $T'_{L'}$. The universal property of scalar restrictions yields a canonical isomorphism ${\rm Hom}_{F'}(\M{G}_m, T') \cong {\rm Hom}_{F}(\M{G}_m, T)$. Tensoring with $K$ over $\M{Q}_p$ gives an isomorphism
\[\alpha: {\rm Hom}_{L'}(\M{G}_m, T'_{L'}) \cong {\rm Hom}_{L}(\M{G}_m, T_L).\]
Note that a priori this is an isomorphism between morphism sets of schemes, but it restricts to a bijection between the sets of group scheme morphisms
\[\alpha: X_*(T'_{L'}) \cong X_*(T_L).\] 
Then $\alpha$ can be extended to a morphism between the rational cocharacter groups.

\lem{}{}{
Let $\psi \in X_*(T'_{L'})$ and $\unif \in L'$ be any uniformizer. Then
\[\alpha(\psi)(\unif) = \phi(\psi(\unif)).\]
}

\prooof
Although this lemma should not come as a surprise, we have to consider the maps on the level of morphism between the underlying algebras, if we want to be precise at this point. Fist of all fix an isomorphism $\MC{O}_{T}(T) \cong F[\{x_i\}_i]/(\{f_j\}_j)$ and a basis $\{\varepsilon_k\}_k$ of $F$ over $F'$. Set $x_i = \sum_k y_{i,k} \varepsilon_k$ (with further variables $y_{i,k}$) and let $f_{j,k}$ be polynomials in the variables $y_{i,k}$ such that $\sum_k f_{j,k}(\{y_{i,k}\}_{i, k}) \varepsilon_k = f_j(\{x_i\}_i)$. Then the explicit description of scalar restrictions gives
\begin{align*}
 \MC{O}_{T'_{L'}}(T'_{L'}) & =  F'[\{y_{i,k}\}_{i,k}]/(\{f_{j,k}\}_{j, k}) \otimes_{F'} L' \\
 & =  L'[\{y_{i,k}\}_{i,k}]/(\{f_{j,k}\}_{j, k}).
\end{align*}
Then consider $\psi: \M{G}_{m, L'} \to T'_{L'}$ and its image of $\unif$:
\begin{align*}
 \psi_0^{\#}:  L'[\{y_{i,k}\}_{i,k}]/(\{f_{j,k}\}_{j, k}) & \;=\; F'[\{y_{i,k}\}_{i,k}]/(\{f_{j,k}\}_{j,k}) \otimes_{F'} L'\\
    & \longrightarrow  F'[t^{\pm 1}] \otimes_{F'} L' \\
    & \;=\;  L'[t^{\pm 1}]  \\[2mm]
   (\{y_{i,k} - \delta_{i,k}\}_{i, k}) & \;=\; (\{y_{i,k} \otimes 1 - 1 \otimes \delta_{i,k}\}_{i, k})  \\
    & \;=\; (\psi_0^{\#})^{-1}(t \otimes 1 - 1 \otimes \unif) \\
    & \;=\; (\psi_0^{\#})^{-1}(t - \unif)
\end{align*}
(for certain elements $\delta_{i, k} \in L'$). Under the given bijection we get for $\alpha(\psi): \M{G}_{m, L} \to T_L$: 
\begin{align*}
 \alpha(\psi)^{\#}: F[\{x_i\}_i]/(\{f_j\}_j) \otimes_{F} L & \longrightarrow F[t^{\pm 1}] \otimes_{F} L \\
    & \;=\;  L[t^{\pm 1}] \\[2mm]
 (\{x_i \otimes 1 - \sum\nolimits_k \varepsilon_k \otimes \delta_{i,k}\}_i) & \;=\; \\
 (\{\sum\nolimits_k y_{i,k}\varepsilon_k \otimes 1 - \varepsilon_k \otimes \delta_{i,k}\}_i) & \;=\; (\alpha(\psi)^{\#})^{-1} (t \otimes 1 - 1 \otimes \unif) \\
   & \;=\; (\alpha(\psi)^{\#})^{-1}(t - \unif)
\end{align*}
But the point $(\{y_{i,k} - \delta_{i,k}\}_{i, k}) \in T'(L') \subset \left({\rm Res}_{F/F'}(G)\right)(L')$ is exactly mapped to $(\{x_i \otimes 1 - \sum_k \varepsilon_k \otimes \delta_{i,k}\}_i) \in T(L) \subset G(L)$ under the map $\phi: \left({\rm Res}_{F/F'}(G)\right)(L') \to G(L)$ (as both points have the same coordinates but the first one has them as a linear combination of the $\varepsilon_k$). Thus
\[\alpha(\psi)(\unif) = \phi(\psi(\unif)),\]
where we identify as usual the element $\unif \in L' \subset L$ with the ideal $(t - \unif)$ in $\Sp L'[t^{\pm 1}]$ respectively in $\Sp L[t^{\pm 1}]$.
\exit

\prop{\label{prop:ExplicitRamifiedDiagram}}{}{
Using the notation as above, let $e$ be the ramification index of $F$ over $F'$. Then the following diagrams commute:
\[
\begin{xy}
 \xymatrix{
   \left({\rm Res}_{F/F'}(G)\right)(L') \ar[r]^-{\phi} \ar[d]_{\nu_L} & G(L) \ar[d]^{\nu_L} & \left({\rm Res}_{F/F'}(G)\right)(L') \ar[r]^-{\phi} \ar[d]_{\mu_L} & G(L) \ar[d]^{\mu_L} \\
  X_*(T'_{L'})_{\M{Q}}/ \Omega \ar[r]^{e \alpha} &  X_*(T_L)_{\M{Q}}/ \Omega  & X_*(T'_{L'})_{\M{Q}}/ \Omega \ar[r]^{e \alpha} & X_*(T_L)_{\M{Q}}/ \Omega
   }
\end{xy} 
\]
}

\prooof
Consider first the Newton points: \\
Let $b' \in \left({\rm Res}_{F/F'}(G)\right)(L')$ and choose $c' \in \left({\rm Res}_{F/F'}(G)\right)(L')$, $\nu' \in X_*(T'_{L'})_{\M{Q}}$, $n' \in \M{N}$ and $\unif' \in L'$ such that: \\
i) $n'\nu' \in X_*(T'_{L'})$. \\
ii) ${\rm Int}(c') \circ (n'\nu')$ is defined over the fixed field ${L'}^{\langle\sigma^{n'}\rangle}$ of $\sigma^{n'}$ on $L'$. \\
iii) $c' \cdot (b'\sigma)^{n'} \cdot c'^{-1} = c' \cdot (n'\nu')(\unif') \cdot c'^{-1} \cdot \sigma^{n'}$. \\
Applying $\phi$ to all elements and denoting $b = \phi(b')$, $c'' = \phi(c')$, $n = n'$ one gets
\[c'' \cdot (b\sigma)^n \cdot {c''}^{-1} = c'' \cdot \phi((n'\nu')(\unif')) \cdot {c''}^{-1} \cdot \sigma^n.\]
Using the lemma for $\psi = n'\nu'$ we have
\[c'' \cdot (b\sigma)^n \cdot {c''}^{-1} = c'' \cdot (\alpha(n'\nu'))(\unif') \cdot {c''}^{-1} \cdot \sigma^n = c'' \cdot (n\alpha(\nu'))(\unif') \cdot {c''}^{-1} \cdot \sigma^n.\]
Now there is a uniformizer $\unif \in L$ and an element $\gamma \in L$ such that
\[\unif^e = \gamma \cdot \unif' \cdot \sigma(\gamma)^{-1}.\]
Then conjugating $b$ further with $(\alpha(n'\nu'))(\gamma)$, i.e. setting $c = c'' \cdot (\alpha(n'\nu'))(\gamma)^{-1}$ gives
\[c \cdot (b\sigma)^n \cdot c^{-1} = c \cdot (\alpha(n'\nu'))(\unif^e) \cdot c^{-1} \cdot \sigma^n = c \cdot (ne\alpha(\nu'))(\unif) \cdot c^{-1} \cdot \sigma^n.\]
This checks that the elements $c$, $\nu = e \cdot \alpha(\nu')$, $n$ and $\unif$ fulfill the third condition used to describe the Newton point of $b$. \\
It remains to show that these elements also satisfy the two remaining conditions: \\
The first condition is obvious as $\alpha$ is a map between the integral cocharacter groups and so is $e \cdot \alpha$. Hence $ne\alpha(\nu') = e \alpha(n'\nu') \in X_*(T_L)$. \\
To show the second condition note that ${\rm Int}((\alpha(n'\nu'))(\gamma)) \circ (ne\alpha(\nu')) = ne\alpha(\nu')$ as we conjugate inside the torus $T_L$. As rising to the $e$th power is even defined over $\M{F}_p$ it remains to show that ${\rm Int}(c'') \circ \alpha(n'\nu')$ is defined over the fixed field of $\sigma^n$ on $L$. But the universal property of scalar restrictions respect being defined over the fixed field of $\sigma^n$ on the respective ground field. Hence a cocharacter already defined over ${L'}^{\langle\sigma^n\rangle}$ is mapped via $\alpha$ to a cocharacter defined over $L^{\langle\sigma^n\rangle}$. As conjugating with $c''$ only changes the chosen maximal torus, ${\rm Int}(c'') \circ \alpha(n\nu)$ is defined over $L^{\langle\sigma^n\rangle}$ if ${\rm Int}(c') \circ (n'\nu')$ is defined over ${L'}^{\langle\sigma^n\rangle}$. But this was assumed. \\
This shows that $\nu = e \alpha(\nu')$ is indeed the Newton point of $b = \phi(b')$. \vspace{1mm} \\
The case of Hodge points is easier: \\
Let $b' \in \left({\rm Res}_{F/F'}(G)\right)(L')$ as before and choose $c'_1, c'_2 \in \left({\rm Res}_{F/F'}(G)\right)(O_{L'})$, a representative of its Hodge point $\mu' \in X_*(T'_{L'})_{\M{Q}}$ and a uniformizer $\unif' \in L'$ such that
\[b' = c'_1 \cdot \mu'(\unif') \cdot c'_2.\]
Applying $\phi$ and denoting $b = \phi(b')$, $c_1 = \phi(c'_1)$ and $c_2 = \phi(c_2)$ gives together with the lemma for $\psi = \mu'$:
\[b = c_1 \cdot \phi(\mu'(\unif')) \cdot c_2 = c_1 \cdot \alpha(\mu')(\unif') \cdot c_2.\]
Note that $\phi$ maps $\left({\rm Res}_{F/F'}(G)\right)(O_{L'})$ to $G(O_L)$, hence $c_1$ and $c_2$ lie in the correct group. \\
Pick now any uniformizer $\unif \in L$. Because of $\unif^e \cdot \unif'^{-1} \in O_L$ there is an element $c_0 \in T(O_L)$ such that
\[\alpha(\mu')(\unif') = \alpha(\mu')(\unif^e) \cdot c_0 = e\alpha(\mu')(\unif) \cdot c_0.\]
Hence
\[b = c_1 \cdot e\alpha(\mu')(\unif) \cdot c_0 c_2 \]
and $e\alpha(\mu')$ is the Hodge point of $b = \phi(b')$.
\exit

\thm{\label{thm:ExplicitRamifiedBound}}{}{
Let $F', F$, $L', L$ and $T \subset G$ be as above and fix a Hodge point $\mu'^1 \in X_*(T'_{L'})_{\M{Q}}$. Let $\mu^1 = e\alpha(\mu'^1) \in X_*(T_L)_{\M{Q}}$ be the image of $\mu'^1$. Then choose a constant $C > 0$ such that for any $N > C$ and any $\un{\mu} = (\mu^1, \mu^2, \ldots, \mu^N) \in (X_*(T_L)_{\M{Q}})^N$ (where the first entry is fixed) the set $H_{\un{\mu}, N}(G) \subset G(L)$ lies inside a Newton stratum. \\
Then for any $N > C$ and any $\un{\mu}' = (\mu'^1, \mu'^2, \ldots, \mu'^N) \in (X_*(T'_{L'})_{\M{Q}})^N$ with first entry the fixed Hodge point $\mu'^1$ the set $H_{\un{\mu'}, N}({\rm Res}_{F/F'}(G)) \subset ({\rm Res}_{F/F'}(G))(L')$ lies inside a Newton stratum.
}

\prooof
Fix $N > C$ and $\un{\mu'} = (\mu'^1, \mu'^2, \ldots, \mu'^N)$ and choose any $L'$-valued point $b' \in H_{\un{\mu'}, N}({\rm Res}_{F/F'}(G))$. Then the commutativity of the right diagram in the proposition gives for $b = \phi(b')$ and any $k = 1, \ldots, N$ (using the ramification index $e$ as in the proposition):
\[\mu_L((b\sigma)^k) = \mu_L(\phi((b'\sigma)^k)) = e\alpha(\mu_{L'}((b'\sigma)^k)) = e\alpha(\mu'^k).\]
For $\mu^k = e\alpha(\mu'^k)$ we have thus $b \in H_{(\mu^1, \mu^2, \ldots, \mu^N), N}(G)$. By assumption this set lies inside a Newton stratum. Hence $\un{\mu}'$ determines the Newton point $\nu_L(b)$. But using the commutativity of the left diagram
\[\nu_L(b) = \nu_L(\phi(b')) = e\alpha(\nu_{L'}(b')).\]
As $e\alpha$ is injective, $\nu_{L'}(b')$ depends only on $\un{\mu}'$ but not on the chosen element $b'$.
\exit

\prop{\label{prop:ExplicitRamifiedBoundGL2}}{}{
Let $F/F'$ be as above. Then for any pair of Hodge points $\un{\mu}' = (\mu'^1, \mu'^2)$ of the group ${\rm Res}_{F/F'}(GL_2)$ the stratum $H_{\un{\mu}', 2}({\rm Res}_{F/F'}(GL_2))$ lies in some Newton stratum.
}

\prooof
Use the previous proposition for $C = \frac 32$, $N = 2$, $G = GL_2$. Then our assumptions are met by Proposition \ref{prop:ExplicitGL2}.
\exit

\rem{}{
A slightly weaker version of this proposition was already shown by Andreatta and Goren \cite[thm. 9.2]{AndreattaGoren}. We will explain how their theorem compares to our result in section \ref{subsec:AndreattaGoren}.
}

\subsection{Restriction of scalars for unramified extensions}\label{subsec:ExplicitUnramified}
The goal of this section is similar to the previous one but for scalar restrictions via unramified extensions. \\
Similarly to above fix a tower of fields $\M{Q}_p \subset F' \subset F \subset \ov{\M{Q}}_p$ with $F'$ and $F$ finite over $\M{Q}_p$ but now $F/F'$ unramified of degree $f$. Let as in \S \ref{subsec:Setting} $L = F \cdot K = F' \cdot K$. \\
Let $G$ be a connected reductive group over $F$. Then by definition of scalar restrictions there is an isomorphism of abstract groups
\[\phi: \left({\rm Res}_{F/F'}(G)\right)(L) \cong G(F \otimes_{F'} L) = \prod_{\tau \in {\rm Hom}_{F'}(F, L)} G(L).\]
Fixing a maximal torus $T \subset G$ as in section \ref{subsec:Setting} gives a maximal torus $T' = {\rm Res}_{F/F'}(T) \subset {\rm Res}_{F/F'}(G)$. Then $\phi$ restricts to $\phi: T'(L) \to \prod_{\tau} T(L)$. \\
Note that $\sigma$ cyclically permutes the $\tau$-factors.

\lem{\label{lem:ExplicitPreBoundUnram}}{}{
Fix a Hodge point $\mu^1 \in X_*(T')_{\M{Q}}$, an integer $d > 0$ and an element $\vartheta \in {\rm Hom}_{F'}(F, L)$. Then there is a constant $C(\mu^1, d)$ depending only on $\mu^1$, $d$ and $\vartheta$ with the following property: For any element $b' \in G(L)$ such that there is an element $b \in {\rm Res}_{F/F'}(G)(L)$ with Hodge point $\mu^1$ and $\mu(b') = \mu(((b\sigma)^d)_{\vartheta})$ (denoting the $\vartheta$-component of $(b\sigma)^d$ under $\phi$ by  $((b\sigma)^d)_{\vartheta}$), the Hodge points of $(b'\sigma)^i$ for $i = 1, \ldots, C(\mu^1, d)$ determine the Newton point of $b'$.
}

\prooof
By Lemma \ref{lem:GeneralHodgeBound} there are only finitely many possibilities for $\mu((b\sigma)^d)$. As $\phi$ maps ${\rm Res}_{F/F'}(G)(O_L)$ into $\prod_{\tau} G(O_L)$, $\mu((b\sigma)^d)$ determines $\mu(\phi((b\sigma)^d))$ and in particular $\mu(((b\sigma)^d)_{\vartheta})$. Applying Theorem \ref{thm:GeneralBound} to the Hodge strata defined by each of them and taking the maximum over all appearing constants gives the desired $C(\mu^1, d)$.
\exit

\thm{\label{prop:ExplicitBoundUnramified}}{}{
Let the $C(\mu^1, d)$ be the constants of the previous lemma. \\ 
i) If $G$ is unramified and has a hyperspecial vertex, then any stratum $H_{\un{\mu}, C(\mu^1, f) \cdot f}({\rm Res}_{F/F'}(G))$ where the first entry of $\un{\mu}$ equals $\mu^1$ lies in some Newton stratum. \\
ii) For general $G$ choose an integer $\xi > 0$ such that $\xi \nu(b)$ lies in the integral cocharacter group for each element $b \in {\rm Res}_{F/F'}(G)(L)$. Then any stratum $H_{\un{\mu}, C(\mu^1, \xi f) \cdot \xi f}({\rm Res}_{F/F'}(G))$ where the first entry of $\un{\mu}$ equals $\mu^1$ lies in some Newton stratum.
}

\prooof
i) Let $\xi$ be some integer such that the Newton points of all elements in ${\rm Res}_{F/F'}(G)$ lie in the $\frac 1\xi X_*(T')$. Then by passing to some totally ramified field extension of $L$ of degree $\xi$, we may wlog. assume that the Newton points of all appearing elements lie in $X_*(T')$ (similarly to the argument in Proposition \ref{prop:ExplicitGL2} but now using Proposition \ref{prop:BaseChangeGeneral}). \\
Consider now an element $b \in H_{\un{\mu}, C(\mu^1, f) \cdot f}({\rm Res}_{F/F'}(G))$ and let $b'_{\tau} \coloneqq ((b\sigma)^f)_\tau \in G(L)$. As $\sigma^f$ is an endomorphism of $\prod_{\tau} G(L)$ which fixes each component, we have $((b\sigma)^{fi})_\tau = (b'_\tau\sigma^f)^i$ for each $i > 0$ and each $\tau$. In particular the fixed Hodge points of $b$ determine the Hodge points $\mu(b'_\tau\sigma), \mu((b'_\tau\sigma)^2), \ldots$, $\mu((b'_\tau\sigma)^{C(\mu^1, f)})$. By choice of $C(\mu^1, f)$ this information suffices to determine the Newton point of $b'_{\vartheta} \in G(L)$ for some $\vartheta$. \\
As the Frobenius $\sigma$ cyclically permutes the $\tau$-factors, we may choose $c \in {\rm Res}_{F/F'}(G)(L)$, $\tilde{b} \in G(L)$ with $\tilde{b} = \nu'(\unif)$ for some $\nu' \in X_*(T)$ and some uniformizer $\unif \in L$ such that
\[\phi(cb\sigma(c)^{-1}) =  (\tilde{b})_{\tau}.\]
Write $\phi(b) = (b_\tau)_{\tau}$ and $\phi(c) = (c_{\tau})_{\tau}$. Then $c_{\tau} b_{\tau} \sigma(c_{\sigma\tau}^{-1}) = \tilde{b}$ in $G(L)$. Thus as $f = [F : F']$ equals the number of embeddings $\tau$ (recall $L$ is the maximal unramified extension of $F$)
\[b'_{\tau} = ((b\sigma)^f)_\tau = c_{\tau}^{-1} \cdot (\tilde{b}\sigma)^f \cdot c_{\tau} \in G(L).\]
In particular the Newton point of $b'_{\vartheta}$ determines the Newton point of $(\tilde{b}\sigma)^f$, hence the one of $\tilde{b}$. But by definition of $\tilde{b}$, $\nu(\tilde{b})$ gives directly the Newton point of $b$ itself. But this was our goal. \\
ii) Instead of changing the base field replace $b$ by $(b\sigma)^{\xi}$. Then we find a $\tilde{b}$ with the same properties as in ii) and the same arguments show that the Hodge points $\mu((b\sigma)^{\xi f}), \mu((b\sigma)^{2\xi f}), \ldots, \mu((b\sigma)^{C(\mu^1, \xi f) \cdot \xi f})$ determine $\nu(((b\sigma)^{\xi f})_\vartheta)$ and hence $\nu(b)$.
\exit

\cor{}{}{
Let $F/F'$ be as in the proposition. Then any stratum $H_{\un{\mu}, 2f}({\rm Res}_{F/F'}(GL_2))$ lies in some Newton stratum.  
}

\prooof
Apply part i) of the previous proposition for $G = GL_2$ and note that due to Proposition \ref{prop:ExplicitGL2} we may choose $C(\mu^1, f) = 2$.
\exit

\rem{}{
i) Note that the Newton points are not contained in $\frac 12 X_*(T')$, but in $\frac 1{2f} X_*(T')$. \\
ii) Contrary to the totally ramified situation the Hodge points of $b$ and $(b\sigma)^2$ do not suffice to determine the Newton point. To see this consider ${\rm Res}_{F/\M{Q}_p}(GL_2)$ with $F/\M{Q}_p$ unramified of degree $3$ and $\mu^1 = ((1, 0), (1, 0), (1, 0))$, $\mu^2 = ((1, 1), (2, 0), (2, 0))$. Then one easily checks that the element
\[b = \left(\lmat{cc} \unif & 0 \\ 0 & 1 \rmat, \lmat{cc} 0 & 1 \\ \unif & 0 \rmat, \lmat{cc} 1 & 0 \\ 0 & \unif \rmat\right) \in \prod_{\tau \in {\rm Hom}_{\M{Q}_p}(F, K)} GL_2(K) \cong {\rm Res}_{F/\M{Q}_p}(GL_2)(K)\]
has Newton point $\nu(b) = ((\frac 12, \frac 12), (\frac 12, \frac 12), (\frac 12, \frac 12))$. But the results of section \ref{subsec:GorenOort} imply that the 'generic' element in $H_{(\mu^1, \mu^2), 2}({\rm Res}_{F/\M{Q}_p}(GL_2))$ has Newton point $((\frac 23, \frac 13), (\frac 23, \frac 13), (\frac 23, \frac 13))$. 
}

\subsection{Restriction of scalars for arbitrary finite field extensions}\label{subsec:ExplicitFieldExt}
We will connect the previous two sections \ref{subsec:ExplicitRamified} and \ref{subsec:ExplicitUnramified} to get corresponding statements for arbitrary scalar restrictions. 

\thm{}{}{
Let $\M{Q}_p \subset F' \subset F \subset \ov{\M{Q}}_p$ be finite field extensions inside $\ov{\M{Q}}_p$. Let $G$ be a connected reductive group over $F$ with maximal torus $T$ and fix a Hodge point $\mu'^1 \in X_*(T')_{\M{Q}}$ where $T' = {\rm Res}_{F/F'}(T) \subset {\rm Res}_{F/F'}(G)$. Then one can give an explicit bound $C$ depending only on the constants for $G$ itself and on invariants of the extension $F/F'$, such that for every tuple $\un{\mu}' = (\mu'^1, \ldots, \mu'^N)$ with $N > C$ and first entry the chosen Hodge point $\mu'^1$ the stratum $H_{\un{\mu}', N}({\rm Res}_{F/F'}(G))$ lies in some Newton stratum.
}

\prooof
Let $F''$ be the maximal unramified extension of $F'$ inside $F$ and $f = [F'':F]$. Denote $L' = F' \cdot K = F'' \cdot K$ and $L = F \cdot K$. Then
\[{\rm Res}_{F/F'}(G) = {\rm Res}_{F''/F'}({\rm Res}_{F/F''}(G))\]
and we have the following isomorphisms of groups:
\begin{align*}
 \phi: ({\rm Res}_{F/F'}(G))(L') & \;=\; ({\rm Res}_{F''/F'}({\rm Res}_{F/F''}(G)))(L') \\
 & \stackrel{\phi''}{\longrightarrow} \prod_{\tau \in {\rm Hom}_{F'}(F'', L')}({\rm Res}_{F/F''}(G))(L') \stackrel{\phi'}{\longrightarrow} \prod_{\tau \in {\rm Hom}_{F'}(F'', L')} G(L)
\end{align*}
where $\phi''$ resp. $\phi'$ are the isomorphisms of section \ref{subsec:ExplicitUnramified} resp. section \ref{subsec:ExplicitRamified}. 
Now fix the Hodge point $\mu'^1 \in X_*(T')_{\M{Q}}$, an element $\vartheta \in {\rm Hom}_{F'}(F'', L')$ and $\xi > 0$ such that all Newton points for ${\rm Res}_{F/F'}(G)$ are contained in $\frac 1{\xi}X_*(T')$. As in Lemma \ref{lem:ExplicitPreBoundUnram} we may find a constant $C(\mu'^1, \xi f)$ satisfying the very same properties as stated there, but now for the morphism $\phi$ considered here. Then applying Theorem \ref{thm:ExplicitRamifiedBound} we see that this $C(\mu'^1, \xi f)$ also satisfies the property of Lemma \ref{lem:ExplicitPreBoundUnram} for the group ${\rm Res}_{F/F''}(G)$, i.e. when considering only the morphism $\phi''$. But this means that all assumptions of Theorem \ref{prop:ExplicitBoundUnramified} are met for the unramified extension $F/F''$, and we may take $C = C(\mu'^1, \xi f) \cdot \xi f$. \\
Note that if $G$ is unramified and contains a hyperspecial vertex the same argument works with $\xi = 1$.
\exit

\prop{\label{prop:ExplicitBoundResGL2}}{}{
Let $\M{Q}_p \subset F' \subset F \subset \ov{\M{Q}}_p$ be finite field extensions inside $\ov{\M{Q}}_p$ and let $f$ be the degree of the maximal unramified extension of $F'$ inside $F$ (as in the proof of the previous theorem). Then for any tuple of Hodge points $\un{\mu} = (\mu^1, \mu^2, \ldots, \mu^{2f})$ of ${\rm Res}_{F/F'}(GL_2)$ the stratum $H_{\un{\mu}, 2f}({\rm Res}_{F/F'}(GL_2))$ lies in some Newton stratum.
}

\prooof
Use the previous theorem in the special case $G = GL_2$. Then we may choose $C = 2f$ by Proposition \ref{prop:ExplicitGL2}.
\exit

\subsection{\texorpdfstring{$\mathbf{SL_n}$ for $\mathbf{n \geq 3}$}{SL(n) for n > 2}}\label{subsec:ExplicitSLn}
\prop{}{}{
For each $n \geq 1$ and arbitrary $L$, there are strata $H_{\un{\mu}, n-1} \subset SL_n(L)$ which are not contained in any Newton stratum.
}

\prooof
Consider for a chosen uniformizer $\unif$ the matrices:
\[b_1 = \lmat{cccccc} 0 & 0 & \ldots & 0 & 0 & (-1)^{n-1} \unif^{n-1} \\ \unif^{-1} & 0 & \ldots & 0 & 0 & 0 \\ 0 & \unif^{-1} & \ldots & 0 & 0 & 0 \\ \vdots & \vdots & \ddots & \vdots & \vdots & \vdots \\ 0 & 0 & \ldots & \unif^{-1} & 0 & 0 \\ 0 & 0 & \ldots & 0 & \unif^{-1} & 0 \rmat \]
\[b_2 = \lmat{cccccc} 0 & 0 & \ldots & 0 & (-1)^n \unif^{n-1} & 0 \\ \unif^{-1} & 0 & \ldots & 0 & 0 & 0 \\ 0 & \unif^{-1} & \ldots & 0 & 0 & 0 \\ \vdots & \vdots & \ddots & \vdots & \vdots & \vdots \\ 0 & 0 & \ldots & \unif^{-1} & 0 & 0 \\ 0 & 0 & \ldots & 0 & 0 & \unif^{-1} \rmat \]
Let $T \subset SL_n$ be the diagonal torus. We will give the Hodge and Newton points via giving the image of $\unif$ on some representative in $X_*(T)_{\M{Q}}$ (for $1 \leq i \leq n - 1$): 
\[\mu((b_1\sigma)^i)(\unif) = {\rm diag}\left((\underbrace{\unif^{n-i}, \unif^{n-i}, \ldots, \unif^{n-i}}_{i}, \underbrace{\unif^{-i}, \ldots, \unif^{-i}}_{n-i})\right) = \mu((b_2\sigma)^i)(\unif)\]
\[\mu((b_1\sigma)^n)(\unif) = (1, 1, \ldots, 1) \neq \mu((b_2\sigma)^n)(\unif) = (\unif^n, 1, 1, \ldots, 1, \unif^{-n})\]
\[\nu(b_1)(\unif) = (1, 1, \ldots, 1) \neq \nu(b_2)(\unif) = (\unif^{\frac 1{n-1}}, \unif^{\frac 1{n-1}}, \ldots, \unif^{\frac 1{n-1}}, \unif^{-1})\]
Note that the last two inequality signs hold even after taking the quotient by $\Omega = S_n$.
\exit

\rem{}{
i) The same statement holds for $G = GL_n$ or $G = PGL_n$ using the same matrices. \\
ii) These examples restrict the cases where one can hope that the first two Hodge points already define the Newton point: No connected reductive group containing a subgroup $SL_3$ or $PGL_3$ has this property.
}

\section{Comparison to stratifications defined by Goren-Oort and \texorpdfstring{\\}{} Andreatta-Goren}\label{sec:OtherStrata}
\subsection{Connection to the work of Goren and Oort}\label{subsec:GorenOort}
Fix a totally real extension $\tilde{F}$ over $\M{Q}$ such that the corresponding extension $F$ over $\M{Q}_p$ is unramified of degree $g$. We consider quadruples $(A, \lambda, \iota, \alpha)$ consisting of an abelian variety $A$ over $k$ (which we assume to be algebraically closed), a principal polarization $\lambda$, $\iota: O_{\tilde{F}} \to End(A)$ fixed by the Rosati involution associated to $\lambda$ and a full symplectic level-$n$-structure $\alpha$. The  moduli space of such tuples is representable by a regular irreducible variety $\MC{M}_n$ over $k$. \\
To such quadruples Goren and Oort associate in \cite{GorenOortStrata} a discrete invariant $\tau(A)$, the \emph{type} of $A$ and study the corresponding stratification $W_{\tau}^0$ on $\MC{M}_n$. We will explain here a group theoretic description of the type and give a conceptual rather than computational proof of a weaker version of \cite[thm. 5.4.11]{GorenOortStrata} about the generic Newton point on the stratum defined by a fixed type $\tau$. \\
Throughout this section $\kappa(\cdot)$ always denotes the residue field of a local field and we fix isomorphisms
\[\M{Z}/g\M{Z} \cong {\rm Hom}_{\M{F}_q}(\kappa(F), k) \cong {\rm Hom}_{\M{Z}_p}(O_F, O_K) \cong {\rm Hom}_{\M{Q}_p}(F, K)\]
(where $K = {\rm Frac}(W(k))$ as in section \ref{subsec:Setting}) such that the homomorphism associated to $i+1$ is obtained by composing the homomorphism associated to $i$ with $\sigma$.

\defi{}{(\cite[def. 2.1.1]{GorenOortStrata})}{
Let $(A, \lambda, \iota, \alpha) \in \MC{M}_n(k)$. Then $\kappa(F)$ acts via $\iota$ on $H^0(\Omega^1_A)$ and hence on the $k$-vector space $ker(V: H^0(\Omega^1_A) \to H^0(\Omega^1_A))$ where $V$ denotes the Verschiebung map. This kernel decomposes as a direct sum of subspaces on which the $\kappa(F)$-action is given via a character in ${\rm Hom}_{\M{F}_q}(\kappa(F), k) = \M{Z}/g\M{Z}$. Then the \emph{type} $\tau(A)$ of $(A, \lambda, \iota, \alpha)$ is the subset of $\M{Z}/g\M{Z}$ consisting of all characters for which the corresponding subspace is non-trivial.
}

The Frobenius morphism $F$ on the rational Dieudonn\'e-module corresponding to $(A, \lambda, \iota, \alpha)$ is given by an element (or rather a ${\rm Res}_{F/\M{Q}_p}(GL_2)(O_K)$-$\sigma$-conjugacy class) $b \in {\rm Res}_{F/\M{Q}_p}(GL_2)(K)$. Using as in section \ref{subsec:ExplicitUnramified}
\[{\rm Res}_{F/\M{Q}_p}(GL_2)(K) \cong \prod_{\tau \in {\rm Hom}_{\M{Q}_p}(F, K)} GL_2(K) = \prod_{i \in \M{Z}/g\M{Z}} GL_2(K)\]
we will view $b$ as a tuple $(b_i)_i$ on the right-hand side. \\
Wlog. let $T \subset GL_2$ be the maximal torus of diagonal elements. Then we use again $T' = {\rm Res}_{F/\M{Q}_p}(T)$ as a maximal torus of ${\rm Res}_{F/\M{Q}_p}(GL_2)$. But rather than working with Hodge or Newton points in $X_*(T')_{\M{Q}}$, we use their image under the identification
\begin{align*}
 X_*(T')_{\M{Q}} & = {\rm Hom}_K(\M{G}_{m, K}, T' \times_{\Sp \M{Q}_p} \Sp K)_{\M{Q}} \\
 & = {\rm Hom}_K(\M{G}_{m, K}, \prod_{i \in \M{Z}/g\M{Z}} (T \times_{\Sp F} \Sp K))_{\M{Q}}  = \prod_{i \in \M{Z}/g\M{Z}} X_*(T)_{\M{Q}}
\end{align*}
Then an easy computation gives for $b \in {\rm Res}_{F/\M{Q}_p}(GL_2)(K)$ with components $(b_i)_i \in \prod_{i \in \M{Z}/g\M{Z}} GL_2(K)$ the equality $\mu(b)_i = \mu(b_i)$ under the isomorphism between cocharacter groups above. For Newton points the situation is slightly different as $\sigma$ permutes all $GL_2$-factors. This implies that the image of $\nu(b)$ in each $X_*(T)_{\M{Q}}$-factor is the same. 

\prop{}{}{
Let $(A, \lambda, \iota, \alpha) \in \MC{M}_n(k)$ and let $b \in {\rm Res}_{F/\M{Q}_p}(GL_2)(K)$ represent the action of the Frobenius on the Dieudonn\'e module associated to $A$. Then $i \in \tau(A)$ if and only if $\mu((b\sigma)^2)_i = (1, 1)$.
}
 
\prooof
Let $(D, F, V)$ be the Dieudonn\'e module of $A$. Then there is a canonical isomorphism $H^0(\Omega^1_A) \cong VD/pD$ as $k = \kappa(K)$-module. Furthermore $D$ decomposes as a direct sum of $O_K$-modules $D_i$ on which the $O_F$-action induced by $\iota$ is given by the element $i \in \M{Z}/g\M{Z} = {\rm Hom}_{\M{Z}_p}(O_F, O_K)$ (cf. \cite[\S2.3]{GorenOortStrata}). Each of the $D_i$ is free of rank $2$. As the Verschiebung map $V$ is $\sigma^{-1}$-linear its restriction $V|_{D_i}: D_i \to D$ factors via $V_i: D_i \to D_{i-1}$. Similarly the restriction $F|_{D_i}$ of the Frobenius factors via $F_i: D_i \to D_{i+1}$. Dividing out $p$ we get $2$-dimensional $k$-vector spaces $D_i/pD_i$ and morphisms $\ov{V}_i: D_{i+1}/pD_{i+1} \to D_i/pD_i$, $\ov{F}_i: D_i/pD_i \to D_{i+1}/pD_{i+1}$. By \cite[lemma 2.3.1]{GorenOortStrata} each of the $V_i$ and $F_i$ has cokernel isomorphic to $O_K/pO_K = k$. In particular each $\ov{V}_i$ has a $1$-dimensional kernel and a $1$-dimensional image. \\
Now we can reformulate the condition for $i \in \tau(A)$: By definition this happens if and only if 
\[\ov{V}_{i-1}: \ov{V}_i(D_{i+1}/pD_{i+1}) = V_iD_{i+1}/pD_i \to \ov{V}_{i-1}(D_i/pD_i) = V_{i-1}D_i/pD_{i-1} \hookrightarrow D_{i-1}/pD_{i-1}\]
has non-trivial kernel. As each $\ov{V}_i$ has $1$-dimensional kernel, this happens exactly if $\ov{V}_{i-1} \circ \ov{V}_i: D_{i+1}/pD_{i+1} \to D_{i-1}/pD_{i-1}$ is the zero morphism. Using that the $V_i$ are injective with the stated cokernel, we may rewrite this as
\[V_{i-1} \circ V_i = p \sigma^{-2}: D_{i+1} \to D_{i-1}.\] 
With $V_i \circ F_i = p: D_i \to D_i$ for all $i$ this is equivalent to
\[F_i \circ F_{i-1} = (p^{-1}\sigma^2) \circ V_{i-1} \circ V_i \circ F_i \circ F_{i-1} = p \sigma^2: D_{i-1} \to D_{i+1}\]
But this condition can be reformulated in group theoretic terms: Fixing a suitable $O_K$-basis of $D$ we can view the Frobenius as $F = b\sigma$ for some element $b \in {\rm Res}_{F/\M{Q}_p}(GL_2)(K)$. Then the group decomposition ${\rm Res}_{F/\M{Q}_p}(GL_2)(K) = \prod_{i \in \M{Z}/g\M{Z}} GL_2(K)$ corresponds exactly to the decomposition with respect to the action on the eigenspaces $D_i$, i.e. we have $F_i = b_i\sigma: D_i \to D_{i+1}$. Thus $F_i \circ F_{i-1} = p \sigma^2$ if and only if
\[(b_i\sigma) \cdot (b_{i-1}\sigma) \in GL_2(O_K) \cdot pE_2 \cdot GL_2(O_K)\] 
(denoting the unit matrix by $E_2$) or in other words $\mu(b_i \cdot \sigma(b_{i-1})) = (1, 1)$. But $b_i \cdot \sigma(b_{i-1})$ is nothing else than the $i$-th component of the element $(b\sigma)^2 \in {\rm Res}_{F/\M{Q}_p}(GL_2)(K)$.
\exit

\cor{}{}{
The loci in $\MC{M}_n$ with constant type $\tau$ are exactly the loci where $\mu(b)$ and $\mu((b\sigma)^2)$ are constant, i.e. the variety $W_{\tau}^0$ is the preimage of $H_{\un{\mu}, 2}({\rm Res}_{F/\M{Q}_p}(GL_2))$ for a suitable $\un{\mu}$ under the map
\[\MC{M}_n \to {\rm Res}_{F/\M{Q}_p}(GL_2)(K) / ({\rm change \; of \; basis}) \]
associating to each quadruple the element defining the action of the Frobenius.
}

\prooof
Let $(A, \lambda, \iota, \alpha) \in \MC{M}_n(k)$ and choose a representative $b \in {\rm Res}_{F/\M{Q}_p}(GL_2)$ of the Frobenius morphism on its Dieudonn\'e module. By \cite[lemma 2.3.1]{GorenOortStrata} $\mu(b) \coloneqq ((1, 0), \ldots, (1, 0))$ is constant on all of $\MC{M}_n$, hence does not give any information at all. But as $\mu((b\sigma)^2)_i$ is either $(1, 1)$ or $(2, 0)$ for each $i$, the previous proposition states that the type encodes precisely the same information as $\mu((b\sigma)^2)$.
\exit

To treat the generic Newton point of $W_{\tau} = \bigcup_{\tau' \subset \tau} W_{\tau'}^0 \subset \MC{M}_n$ (c.f. \cite[def. 2.3.5]{GorenOortStrata}), we recall its description found by Goren and Oort:

\defi{}{}{
A subset $\tau' \subseteq \M{Z}/g\M{Z}$ is called spaced if it contains no two consecutive elements. For any $\tau \subset \M{Z}/g\M{Z}$ set
\[\lambda(\tau) = \left\{\begin{array}{ccc} \frac 12 & {\rm if} & g {\rm \; odd \; and \;} \tau = \M{Z}/g\M{Z} \\ 
                                      \frac 1g \max \{|\tau'| \;|\; \tau' \subseteq \tau \, {\rm spaced}\} & {\rm else} & \end{array}\right.\]
\[\beta_{\tau} = \left((1-\lambda(\tau), \lambda(\tau)), \ldots, (1-\lambda(\tau), \lambda(\tau))\right) \in X_*(T')_{\M{Q}} = \prod_{i \in \M{Z}/g\M{Z}} X_*(T)_{\M{Q}}\] 
}

\rem{}{
Note that the given definitions of $\lambda(\tau)$ and $\beta_\tau$ differ slightly from the definitions $5.2.3$ and $5.4.8$ in \cite{GorenOortStrata} as we wish them to be compatible with the notation used in the following section. 
}

\prop{}{}{
Let $\tau \subset \M{Z}/g\M{Z}$ be a type. Then the Newton point $\nu$ of every geometric point in $W_{\tau}$ satisfies $\nu \prec \beta_{\tau}$.
}

\prooof
Note first that for $\tau' \subseteq \tau$ we have $\beta_{\tau} \prec \beta_{\tau'}$. As $W_{\tau} = \bigcup_{\tau' \subseteq \tau} W_{\tau'}^0$ we are reduced to consider only points with type $\tau$. \\
Thus fix any point in $W_{\tau}^0$ and let $b \in {\rm Res}_{F/\M{Q}_p}(GL_2)$ represent the action of the Frobenius morphism on the Dieudonn\'e module. If $g$ is odd and $\tau = \M{Z}/g\M{Z}$ then we have $\mu((b\sigma)^2) = ((1, 1), \ldots, (1, 1))$. Therefore by Mazur's inequality for unramified groups (cf. \cite[thm. 4.2]{RapoRich})
\[\nu(b) = \frac 12 \nu((b\sigma)^2) \prec \frac 12 \mu((b\sigma)^2) = \beta_{\tau}.\]
Consider now the remaining cases and recall
\[\mu(b)_i = (1, 0) \quad {\rm and} \quad \mu((b\sigma)^2)_i = \left\{\begin{array}{ccc} (1, 1) & {\rm if} & i \in \tau \\ (2, 0) & {\rm if} & i \notin \tau \end{array}\right.\]
Fix a spaced subset $\tau' \subseteq \tau$ of maximal cardinality. Replacing $b$ by $\sigma^k(b)$ for suitable $k$ (which does not change the Newton point), we may assume that $0 \in \tau'$ if $\tau'$ is non-empty. Now consider the partition $\MC{P}_{\tau'}$ of $\M{Z}/g\M{Z}$ defined by \\
i) $\{i\} \in \MC{P}_{\tau'}$ if and only if $i \notin \tau'$ and $i+1 \notin \tau'$. \\
ii) $\{i-1, i\} \in \MC{P}_{\tau'}$ if and only if $i \in \tau'$. \\
With Lemma \ref{lem:GLHodgeEstimate} we get
\[\mu((b\sigma)^g)_0 = \mu\left(\prod_{I \in \MC{P}_{\tau'}} \left(\prod_{i \in I} \sigma^{g-i}(b_i)\right)\right) \prec \bigoplus_{I \in \MC{P}_{\tau'}} \mu\left(\prod_{i \in I} \sigma^{g-i}(b_i))\right)\]
(for a suitable enumeration of the elements in $\MC{P}_{\tau'}$ in the expression in the middle). By our choice of the partition we have 
\[\mu\left(\prod_{i \in I} \sigma^{g-i}(b_i)\right) = \left\{\begin{array}{ccc} (1, 0) & {\rm if} & |I| = 1 \\ (1, 1) & {\rm if} & |I| = 2 \end{array}\right. \]
and obtain
\[\mu((b\sigma)^g)_0 \prec (1, 0)^{\oplus g - 2|\tau'|} \oplus (1, 1)^{\oplus |\tau'|} = (g - |\tau'|, |\tau'|).\]
Thus for any $i \in \M{Z}/g\M{Z}$ and any $n > 0$:
\begin{align*}
 \mu((b\sigma)^{ng})_i & = \mu((b\sigma)^i \cdot ((b\sigma)^g)^{n-1} \cdot (b\sigma)^{g-i})_i \\
  &  \prec \bigoplus_{j = i, i-1, \ldots, 1}\mu(b)_j \oplus \mu((b\sigma)^g)_0^{\oplus n-1} \oplus \bigoplus_{j = g, g-1, \ldots, i+1} \mu(b)_j \\
  &  = (1, 0)^{\oplus g} \oplus (g - |\tau'|, |\tau'|)^{\oplus n-1} \\
  &  = (ng - (n-1)|\tau'|, (n-1)|\tau'|).
\end{align*}
Hence we have again by Mazur's inequality
\[\nu(b) \prec \lim_{n \to \infty} \frac 1{gn} \mu((b\sigma)^{gn}) = \frac 1g((g - |\tau'|, |\tau'|), \ldots, (g - |\tau'|, |\tau'|)) = \beta_{\tau}.\]
\exit

\rem{}{
This proposition is a weaker form of \cite[thm. 5.4.11]{GorenOortStrata}. There an explicit calculation using the matrix 
\[A_{\tau} = \left(\lmat{cc} a_i & -\unif \\ 1 & 0\rmat\right)_i \qquad {\rm with} \quad a_i = \left\{\begin{array}{ccc} 0 & {\rm if} & i \in \tau \\ \op{ generic \; and \; non-zero}& {\rm if} & i \notin \tau \end{array}\right.\]
shows that on each irreducible component of $W_{\tau}$ the Newton point $\beta_{\tau}$ indeed occurs and is the generic one. Nevertheless we included the proof above as it gives a more conceptual explanation why one should expect this behavior: The fixed Hodge points $\mu(b)$ and $\mu((b\sigma)^2)$ give by using Lemma \ref{lem:GeneralHodgeBound} upper bounds on $\frac 1n\mu((b\sigma)^n)$ (for $n > 0$). Then we expect the generic Newton point to be the biggest Newton point satisfying Mazur's inequality for all these upper bounds.
}

\subsection{Connection to the work of Andreatta and Goren}\label{subsec:AndreattaGoren}
Fix a totally real extension $\tilde{F}$ over $\M{Q}$ such that the corresponding extension $F$ over $\M{Q}_p$ is totally ramified of degree $g$. We consider quadruples $(A, \lambda, \iota, \varepsilon)$ consisting of an abelian variety $A$ over a field of characteristic $p$, a polarization $\lambda$, $\iota: O_{\tilde{F}} \to End(A)$ and a level-$N$-structure $\varepsilon$ (with $N \geq 4$) satisfying the Deligne-Pappas condition. For a precise treatment see \cite[\S 2.1]{AndreattaGoren}. The moduli space of such tuples is representable by a smooth irreducible variety $\MC{M}(\M{F}_p, \mu_N)$ over $\M{F}_p$. \\
We recall the definition of the two invariants $j$ and $n$ defined in \cite{AndreattaGoren} using the display associated the tuple $(A, \lambda, \iota, \varepsilon)$ defined over an algebraically closed field $k$. For this set as usual $K = {\rm Frac}(W(k))$ and $L = F \cdot K$. Then $L/K$ is a totally ramified extension of degree $g$ and its ring of integers can be identified with $O_L = O_{\tilde{F}} \otimes_{\M{Q}} W(k) \cong W(k)[T]/(h(T))$ for some Eisenstein polynomial $h(T)$ of degree $g$ (cf. \cite[\S4.4]{AndreattaGoren}).
 
\defi{}{(cf. \cite[def. 4.1]{AndreattaGoren})}{
Let $R$ be some $\M{F}_p$-algebra and $W(R)$ the Witt vectors over $R$ with the Frobenius morphism $\sigma$. An $O_{\tilde{F}}$-display (over $R$) is a quadruple $(P, Q, V^{-1}, F)$ with: \\
i) $P$ a projective $O_{\tilde{F}} \otimes W(R)$-module of rank $2$. \\
ii) $Q \subseteq P$ a finitely generated $O_{\tilde{F}} \otimes W(R)$-submodule. \\
iii) $F: P \to P$ linear with respect to $O_{\tilde{F}}$ and $\sigma$-linear with respect to $W(R)$. \\
iv) $V^{-1}: Q \to P$ linear with respect to $O_{\tilde{F}}$ and $\sigma$-linear with respect to $W(R)$. \\
satisfying several conditions as stated in \cite[def. 1]{ZinkDisplay}.
}

\prop{\label{prop:DisplayNormalForm}}{(\cite[prop. 4.10]{AndreattaGoren})}{
Let $(P, Q, V^{-1}, F)$ be an $O_{\tilde{F}}$-display over $k$. Let $\ov{P} = P \otimes_{W(k)} k$ and $\ov{Q}$ the image of $Q$ under the projection $P \to \ov{P}$. Let $\ov{F}: \ov{P} \to \ov{P}$ be the reduction of $F: P \to P$. Then there are $\alpha, \beta \in P$ and uniquely determined $i, j, m \in \M{Z}$ such that \\
i) $P = O_L\alpha \oplus O_L\beta$. \\
ii) The Hodge filtration $\ov{Q} = ker(\ov{F}) \subset \ov{P}$ is given by  
\[\ov{Q} = (T^iO_L/p)\alpha \oplus (T^jO_L/p)\beta \subset \ov{P} = (O_L/p)\alpha \oplus (O_L/p)\beta \]
with $i + j = g$, $0 \leq j \leq i \leq g$. \\
iii) There exists a $c \in O_L^\times$ such that
\[F(\alpha) = T^m\alpha + T^j\beta \qquad , \qquad F(\beta) = c \cdot T^i \alpha\]
with $m \geq j$.
}

\defi{}{}{
For an $O_L$-display $(P, Q, V^{-1}, F)$ define the following two integers: \\
i) $j$ as in the previous proposition. \\
ii) $n = \min\{m, i\}$, with $m$ and $i$ as in the previous proposition. \\
iii) A display $(P, Q, V^{-1}, F)$ is of type $(j, n)$ if these integers are the constants defined in i) and ii). \\
Furthermore let $\lambda(n) = \min\{\frac ng, \frac 12\}$.
}

Note that the display associated to a tuple $(A, \lambda, \iota, \varepsilon) \in \MC{M}(\M{F}_p, \mu_N)(k)$ is exactly an $O_{\tilde{F}}$-display over $k$. Thus it makes sense to define $\MC{M}_{j, n} \subset \MC{M}(\M{F}_p, \mu_N)$ as the locus where the display associated to the abelian variety has type $(j, n)$.

\thm{\label{thm:AGNewton}}{(\cite[thm. 9.2]{AndreattaGoren})}{
For any geometric point $(A, \lambda, \iota, \varepsilon) \in \MC{M}_{j, n}(k)$ the slopes of the Newton polygon of $A$ are $\lambda(n)$ and $1 - \lambda(n)$.
}

We will give an alternative proof of this theorem and explain how the type $(j, n)$ is related to certain Hodge points. To do so we will denote the elements of ${\rm Res}_{F/\M{Q}_p}(GL_2)(K)$ as matrices via the isomorphism of groups ${\rm Res}_{F/\M{Q}_p}(GL_2)(K) \cong GL_2(L)$.

\prop{}{}{
Let $(P, Q, V^{-1}, F)$ be an $O_{\tilde{F}}$-display and view $F \in {\rm Res}_{F/\M{Q}_p}(GL_2)(K) = GL_2(L)$ by choice of some basis. Then \\
i) $j$ is the smaller Hodge slope of $F$ (as an element in $GL_2(L)$). \\
ii) $n + j$ is the smaller Hodge slope of $(F\sigma)^2$ (as an element in $GL_2(L)$). \\ 
In particular $\MC{M}_{j, n}$ is exactly the locus where the Frobenius of the display associated to the abelian variety lies in $H_{(\frac 1g(g-j, j), \frac 1g(2g - n - j, n + j)), 2}({\rm Res}_{F/\M{Q}_p}(GL_2))$.
}

\prooof
By Proposition \ref{prop:DisplayNormalForm} we may assume that the display is given in its normal form. In particular one has:
\[F = \lmat{cc} T^m & cT^i \\ T^j & 0\rmat\]
i) Then:
\[ \lmat{cc} 0 & 1 \\ c^{-1} & -c^{-1}T^{m-j}\rmat \cdot F  = \lmat{cc} T^j & 0 \\ 0 & T^i \rmat\]
and the Hodge slopes of $F$ are $j$ and $i$. \\
ii) One computes
\[(F\sigma)^2 = \lmat{cc} c T^{i+j} + T^{2m} & \sigma(c) T^{i+m} \\ T^{j+m} & \sigma(c)T^{i+j} \rmat\]
\textbf{Case 1: $m \geq i$, i.e. $n = i$} \\
Then the following equation yields the Hodge slopes: 
\[\lmat{cc} 1 & -T^{m-j} \\ 0 & 1\rmat \cdot (F\sigma)^2 \cdot \lmat{cc} 1 & 0 \\ -\sigma(c)^{-1}T^{m-i} & 1\rmat = \lmat{cc} c T^{j + i} & 0 \\ 0 & \sigma(c) T^{j + i} \rmat\]
Thus both Hodge slopes of $(F\sigma)^2$ are $i + j = n + j$. \\
\textbf{Case 2: $m < i$} \\
Then the following equation yields the Hodge slopes: 
\[\lmat{cc} 0 & 1 \\ 1 & -cT^{i-m} - T^{m-j} \rmat \cdot (F\sigma)^2 \cdot \lmat{cc} 1 & - \sigma(c) T^{i-m} \\ 0 & 1\rmat = \lmat{cc} T^{j + m} & 0 \\ 0 & -c \sigma(c) T^{2i + j - m} \rmat\]
Thus the smaller Hodge slope of $(F\sigma)^2$ is $m + j = n + j$. \\
The last assertion follows directly from this description of $j$ and $n$. That the slopes get divided by $g$ is due to the fact that we now consider $F \in {\rm Res}_{F/\M{Q}_p}(GL_2)(K)$ and use Proposition \ref{prop:ExplicitRamifiedDiagram} to compare the Hodge points.
\exit

\proof (of theorem \ref{thm:AGNewton}) $\left. \right.$ \\
By Proposition \ref{prop:ExplicitBoundResGL2} all elements in $H_{(\frac 1g(g-j, j), \frac 1g(2g - n - j, n + j)), 2}({\rm Res}_{F/\M{Q}_p}(GL_2))$ have the same Newton point and hence the same holds for $\MC{M}_{j, n}$. To actually compute it use the procedure explained at the end of the proof of Proposition \ref{prop:ExplicitGL2}.
\exit

\bibliographystyle{alpha}
\bibliography{RefinedHodge}

\end{document}